%
%

\magnification=1150

\font\titfont=cmr10 at 12 pt

\font\headfont=cmr10 at 12 pt

\overfullrule=0in

\def\boxit#1{\hbox{\vrule
 \vtop{%
  \vbox{\hrule\kern 2pt %
     \hbox{\kern 2pt #1\kern 2pt}}%
   \kern 2pt \hrule }%
  \vrule}}

  \def\harr#1#2{\ \smash{\mathop{\hbox to .3in{\rightarrowfill}}\limits^{\scriptstyle#1}_{\scriptstyle#2}}\ }

\def\SHF{16}

 \def\GG{{{\bf G} \!\!\!\! {\rm l}}\ }

\def\GL{{\rm GL}}

\def\bll{I \!\! L}

\def\bra#1#2{\langle #1, #2\rangle}
\def\bbf{{\bf F}}
\def\bbj{{\bf J}}
\def\Jtn{{\bbj}^2_n}

\def\ss{\subset}

\def\half{\hbox{${1\over 2}$}}

\def\oa#1{\overrightarrow #1}

\def\dist{{\rm dist}}

\def\log{{\rm log}}
\def\Hess{{\rm Hess}}

\def\tr{{\rm tr}}
\def\max{{\rm max}}

\def\Hom{{\rm Hom\,}}
\def\det{{\rm det}}

\def\Sym{{\rm Sym}^2}

\def\Id{{\rm Id}}

\def\arr{\longrightarrow}

\def\rn{\bbr^n}

\def\Int{{\rm Int}}

\def\Symn{{\Sym(\rn)}}

 \def\cd{{\cal C}}

\def\Theorem#1{\medskip\noindent {\bf THEOREM \bf #1.}}
\def\Prop#1{\medskip\noindent {\bf Proposition #1.}}
\def\Cor#1{\medskip\noindent {\bf Corollary #1.}}
\def\Lemma#1{\medskip\noindent {\bf Lemma #1.}}
\def\Remark#1{\medskip\noindent {\bf Remark #1.}}
\def\Note#1{\medskip\noindent {\bf Note #1.}}
\def\Def#1{\medskip\noindent {\bf Definition #1.}}

\def\Ex#1{\medskip\noindent {\bf Example \bf    #1.}}

\def\pf{\medskip\noindent {\bf Proof.}\ }
\def\qed{\hfill  $\vrule width5pt height5pt depth0pt$}

\def\n{\nabla}

   \def\cp{{\cal P}}

   \def\cn{{\cal N}}
\def\cd{{\cal D}}

\def\cp{{\cal P}}
\def\cf{{\cal F}}

\def\vf{\varphi}

\def\wt{\widetilde}

\def\and{\qquad {\rm and} \qquad}
\def\arr{\longrightarrow}

\def\bbr{{\bf R}}\def\bbh{{\bf H}}
\def\bbc{{\bf C}}

\def\bbp{{\bf P}}

\def\bbm{{\bf M}}
\def\bbg{{\bf G}}

\def\a{\alpha}
\def\b{\beta}
\def\d{\delta}
\def\e{\epsilon}

\def\g{\gamma}

\def\l{\lambda}
\def\o{\omega}

\def\s{\sigma}
\def\x{\xi}

\def\D{\Delta}
\def\L{\Lambda}
\def\bL{{\bf \Lambda}}
\def\G{\Gamma}
\def\O{\Omega}

\def\psh{plurisubharmonic }

\def\lloc{L^1_{\rm loc}}

\def\bo{\partial \Omega}

\def\Symn{\Sym(\rn)}
 
\def\USC{{\rm USC}}
\def\fa{{\rm\ \  for\ all\ }}

\def\cpt{\wt{\cp}}
\def\ft{\wt F}
\def\ob{\overline{\O}}

\def\Fa{{\oa F}}

\def\AA{1}
\def\BB{2}
\def\CC{3}
\def\DD{4}
\def\EE{5}
\def\FF{6}
\def\GGG{7}
\def\HH{8}
\def\II{9}
\def\JJ{10}
\def\KK{11}

\centerline{\titfont  EXISTENCE, UNIQUENESS AND REMOVABLE SINGULARITIES }
\medskip

\centerline{\titfont   FOR NONLINEAR PARTIAL DIFFERENTIAL EQUATIONS }
  
\medskip

\centerline{\titfont   IN GEOMETRY }
  
\bigskip

\centerline{\titfont F. Reese Harvey and H. Blaine Lawson, Jr.$^*$}
\vglue .9cm
\smallbreak\footnote{}{ $ {} \sp{ *}{\rm Partially}$  supported by
the N.S.F. } 

\vskip .2in

\centerline{\bf ABSTRACT} \medskip
  \font\abstractfont=cmr10 at 10 pt

  {{\parindent= .93in

\narrower\abstractfont \noindent
This paper surveys  some recent results
on existence, uniqueness and removable singularities
for fully nonlinear differential equations on manifolds.
The discussion also treats restriction theorems and 
the strong Bellman principle.

}}

\vskip.5in

\centerline{\bf TABLE OF CONTENTS} \bigskip

{{\parindent= .1in\narrower\abstractfont \noindent

\qquad \AA.     Introduction.    \smallskip

\qquad \BB.    Subequations -- a Geometric Approach.    \smallskip

\qquad \CC.    Jet Equivalence of Subequations.   \smallskip

\qquad \DD.  Monotonicity.  \smallskip

\qquad \EE.  Comparison and Strict Approximation.  \smallskip

\qquad \FF.  Removable Singularities.  \smallskip

\qquad \GGG.  Boundary Convexity.  \smallskip

\qquad \HH. The Dirichlet Problem.  \smallskip

\qquad \II.   Restriction Theorems.

\smallskip

\qquad \JJ.   Convex Subequations and the Strong Bellman Principle.

\smallskip

\qquad \KK.  Applications to Almost Complex Manifolds.

}}

{{\parindent= .3in\narrower

\vskip .2in

\hskip .5in 
Appendix A.\    A Pocket Dictionary. 
\smallskip

\hskip .5in 
Appendix B.\    Examples of Basic Monotonicity Cones. 
\smallskip

}}

\vfill\eject


\noindent{\headfont \AA.\  Introduction}
\medskip
 
 Calibrated geometries  are considered generalizations of K\"ahler geometry.
They  resemble K\"ahler geometry in having large families of distinguished
 subvarieties determined by a fixed differential form.  On the other hand,  they seemed at first
 to be unlike K\"ahler geometry in having no
 suitable   analogue of  holomorphic   functions.
 However, it was realized several years ago that the analogues of pluri{\sl sub}harmonic
 functions do exist (in abundance) on any calibrated manifold, and a 
 potential theory was developed in this context  [HL$_{2,3}$].   This led us naturally to the study
 of ``maximal'' or ``extremal'' functions, the analogues of solutions to the 
 homogeneous complex Monge-Amp\`ere equation, first considered by 
 Bremermann  [B] and Walsh [W] and later developed, in the inhomogeneous case,
  by Bedford-Taylor [BT$_*$] and others.
 The techniques and results developed in our study turned out to have
 substantial applications outside of calibrated geometry -- in particular to 
 many of the  highly degenerate elliptic equations which
 appear naturally in  geometry.

 This paper is a survey of those techniques and results.
 We will address questions of existence and uniqueness for the Dirichlet
 Problem, the question of removable singularities for solutions and 
 subsolutions, and the problem of restriction.  The techniques apply
broadly to fully nonlinear (second-order) equations in geometry, and
 in particular, to those which arise ``universally'' on riemannian, hermitian,
 or calibrated manifolds. A  number of examples and applications
 will be discussed, including a proof of the Pali Conjecture on almost complex
 manifolds.  Many more examples appear in the references.
 
 It is conventional in discussing nonlinear differential equations
 to introduce  the notions of a subsolution and supersolution,
 and define a solution to be a function which is both.
 In this paper we adopt an intrinsic approach 
  by specifying a {\sl subset}  $F$ of constraints on the value
 of a function and its derivatives.  
 The classical  subsolutions
 are defined to be the $C^2$-functions $u$ whose 2-jet $(u, Du, D^2u)$ lies in $F$
 at each point.  The set $F$ will be called a {\bf subequation}, and the 
 functions $u$ with $(u, Du, D^2u) \in F$ are called $F$-{\bf subharmonic}.

 The notion of supersolution is captured by the {\bf dual} subequation
 $$
 \ft \  \equiv\  -\{\sim\Int F \} \ =\  \sim\{-\Int F \},
 $$
 and classical solutions $u$ are just those where  $u$ is  $F$-subharmonic and  $-u$ is $\ft$-subharmonic.
 They have the property that $(u, Du, D^2u) \in \partial F$ at each point, since $\partial F = F\cap (\sim \ft)$,
 and they will be called $F$-{\bf harmonic} functions.

 The simplest example is the Laplace equation, where $F=\{ \tr (D^2 u)\geq0\} = \ft$.
 
 The most basic example is the  Monge-Amp\`ere subequation  
 $\cp = \{ D^2u\geq 0\}$ with $\partial \cp \ss \{\det D^2 u=0\}$.
 The dual $\cpt$-subharmonics are the {\sl subaffine functions} (see \BB.1.8).

 Adopting this point of view brings out an  internal duality:
 $$
 \wt{\ft} \ =\ F,
 $$
and enables the roles of $F$ and $\ft$ to be
interchanged  in the analysis.  This  symmetry  is often
 enlightening.  It is particularly  so when discussing the boundary
 geometry necessary for solving the Dirichlet problem.

 A dictionary relating this approach to the more classical one
is given in Appendix A.

The first   step in our  analysis is to extend the notion of $F$-subharmonicity to 
general upper semi-continuous   $[-\infty, \infty)$-valued functions.
This is done in \S 2 where it is noted that these generalized
$F$-subharmonic functions enjoy essentially all the useful properties
of classical subharmonic functions. However,  for this to be meaningful, $F$ must 
 satisfy a  certain  {\sl positivity condition}, corresponding to weak ellipticity.  
 We also require a {\sl negativity condition}, corresponding to weak ``properness''.  
 
For the sake of clarity our exposition will often jump between the two extreme
cases:
\smallskip

(1) \ Constant coefficient (parallel) subequations in $\rn$, and

\smallskip

(2) \ General subequations on manifolds.
 
 \medskip

 In fact, for many equations of interest in geometry and, in particular,
  those which are the principal focus of this survey, 
  these two cases are directly related by the 
 notion of {\bf jet-equivalence}, introduced in \S 3.  This  basic concept plays a fundamental
 role in our work. Jet-equivalence is a certain transformation of  all the variables.
 It can often be quite radical -- turning  mild equations into nasty ones, homogeneous
 equations into inhomogeneous ones, etc.  
 
 As stated, many important nonlinear
 equations on manifolds are locally jet-equivalent, in local coordinates, to
 constant coefficient equations.  In this case  the  results of Slodkowski [S$_1$] and 
 Jensen [J$_1$],  and methods of viscosity theory [CIL], [C] can be applied to prove
{\sl local weak comparison}, and therefore global {\sl weak comparison} ---
the first main step in the analysis of the Dirichlet Problem.

   This leads to another concept of basic importance here: that of a {\bf monotonicity cone},
   introduced in \S 4.
 It gives the approximation tools needed to promote weak comparison to {\sl full comparison}
(see Definition \EE.1)  which,  together with appropriate boundary geometry, yields both  uniqueness
  and existence for the Dirichlet Problem.  A subequation $M$ is called a {\sl monotonicity cone} for a 
  subequation $F$ if 
  $$
  F+M \  \ss\ F
    \eqno{(\AA.1.1)}
    $$
 and each fibre $M_x$,  for $x\in X$, is a convex cone with vertex at the origin. One has that
 $$
  F+M \  \ss\ F\qquad\iff\qquad  \ft+M \  \ss\ \ft,
   $$
   so a monotonicity cone for $F$ is also one for $\ft$.

 Monotonicity cones  play a role in the theory of removable singularities.
For  $M$  as above, we define a
 closed subset $E\ss X$   to be $M$-{\sl polar} if $E = \{x : \psi(x)=\-\infty\}$
 for some $M$-subharmonic function which is smooth on $X-E$.  
 \medskip
 
 \centerline{ \sl
 If $M$ is a monotonicity
 cone for a subequation $F$,}
 
 \centerline{\sl  then
 $M$-polar sets are removable for $F$-subharmonic
 and $F$-harmonic functions on $X$.}

 \medskip\noindent (See Theorems \FF.2.1 and \FF.2.2.)
 This applies, for example, to all branches of the complex Monge-Amp\`ere equation (see \BB.1.10).
Moreover, if a constant pure second-order subequation $F$ in $\rn$ is $M$-monotone, where 
$M \equiv \cp(p) \ss \Symn$ is defined in terms of the ordered eigenvalues by
$
\l_1(A) + \cdots +\l_{[p]}(A) +(p-[p])\l_{p+1}(A) \geq0,
$
then

\medskip

\centerline{
\sl
any closed subset of 
locally finite Hausdorff $p-2$ measure is removable for $F$ and $\ft$.} 

\medskip\noindent
 This applies to the calibration
case.  It generalizes certain results in [CLN], [AGV] and [La$_*$].

Monotonicity cones also play a key role in comparison.  The monotonicity condition (\AA.1.1) is equivalent to
$$
F+ \ft\ \ss\ \wt M.
$$
For  many basic monotonicity cones, the  $\wt M$-subharmonic functions satisfy the 
Zero Maximum Principle (see Appendix B). In such cases, comparison (see \EE.1) comes down
to an {\sl addition theorem}: if $u$ is $F$-subharmonic and $v$ is $\ft$-subharmonic, then 
$u+v$ is $\wt M$ subharmonic.

 There is a last ingredient needed for the Dirichlet Problem -- the necessary boundary
 geometry.  Associated to each subequation $F$, there is a notion of strict $F$-{\bf convexity}
 for oriented hypersurfaces. It is defined in terms of the asymptotic geometry of $F$
 at infinity (see \S 7). It is quite often easy to compute, and it can be expressed directly
 in terms of the second fundamental form.  There are certain equations, like the k-Laplacian
 for $1 < k \leq \infty$ (see \GGG.4(a)), for which all hypersurfaces are strictly $F$-convex.

 A basic result is that:
 \medskip
 \centerline{ \sl If comparison holds for a subequation $F$ on a manifold $X$,}
 
 \centerline{\sl 
 then the Dirichlet Problem is uniquely solvable for $F$-harmonic functions}
 
 \centerline{\sl on every domain $\O  \ss X$ with smooth boundary which is strictly $F$ and $\ft$ convex.}
 \medskip

 Unique solvability for the Dirichlet Problem means that for every $\vf\in C(\bo)$, there exists a unique    $u\in C(\ob)$ such that 
 $$
 u\bigr|_{\O}\in F(\O)
 \and u\bigr|_{\bo} = \vf
 $$
  
  \noindent
  This theorem combines with results discussed above to prove the following
  general result.

 \Theorem{\HH.1.2} {\sl
Let $F$ be a subequation with monotonicity cone $M$.  Suppose that:
\medskip

(i)\ \   $F$ is locally affinely  jet-equivalent to a constant coefficient subequation, and

\medskip

(ii)\ \ $X$ carries a smooth strictly $M$-subharmonic function.
\medskip\noindent
Then existence and uniqueness
hold for the Dirichlet problem for $F$-harmonic functions on any domain
$\O\ss\ss X$  whose boundary
is both strictly $F$-  and   $\ft$-convex.
}
\medskip

The global condition (ii) is essential for a result of this generality. For example, suppose $X$ is a 
riemannian manifold and $F \equiv \{\Hess \,u\geq0\}$,
 where $\Hess\,u$ is the riemannian hessian. Given a  domain $\O\ss \ss X$ 
 with strictly convex boundary, one can  completely change the geometry and topology 
 in the interior of $\O$  without affecting the boundary. The subequation $F$
 continues to satisfy (i), but solutions to the Dirichlet Problem won't exist unless
 (ii) is satisfied.  Another good example is the complex analogue $F=\cp^\bbc$
 on an almost complex hermitian manifold (the homogeneous complex Monge-Amp\`ere
 equation).  Here condition (ii) amounts to the hypothesis that $X$ carries at least
 one strictly plurisubharmonic function.

In homogeneous spaces one can apply a trick of Walsh [W] to establish 
existence without uniqueness.

\Theorem{\HH.1.3} {\sl
Let $X=G/H$ be a riemannian homogeneous space and suppose that $F\ss J^2(X)$ is a 
subequation which is invariant under the natural action of $G$ on $J^2(X)$.
Let $\O\ss\ss X$ be a connected domain whose boundary is both $F$ and $\ft$
strictly convex. Then existence holds for the Dirichlet problem for $F$-harmonic functions on $\O$.
}

\medskip

These results apply to a wide spectrum of equations.  Many examples have been discussed
in [HL$_{4,6,7}$]  and are summarized in \S 2 below.

\smallskip
\noindent
$\bullet$  \ {\bf (Constant Coefficients).}  Theorem \HH.1.3 establishes existence 
for any constant coefficient subequation $F$
in $\rn$, and uniqueness also follows, by \HH.1.2, whenever $F$ has monotonicity cone $M$ and
 there exists a strictly
$M$-subharmonic function on $\ob$.  If $F$ is  pure second-order, for example,
the function $ |x|^2$  works for any $M$, and so uniqueness always  holds. 

For invariant equations on a
sphere, existence always holds by Theorem \HH.1.3. However,  for domains
which do not lie in a hemisphere, where there exists a convex function,
 comparison and its consequences can fail, even for pure second-order equations
(see Appendix D in [HL$_6$]).

\smallskip
\noindent
$\bullet$  \ {\bf (Branches).}  The homogeneous Monge-Amp\`ere equations 
over $\bbr, \bbc$ or $\bbh$ each have branches defined by $\l_k(D^2u) =0$
where $\l_1\leq  \cdots  \leq \l_n$ are the ordered eigenvalues. (See \BB.1.3 and \BB.1.10.)
In fact the equation given by the $\ell^{\rm th}$ elementary symmetric function
$\s_{\ell}(D^2u)=0$ also has $\ell$ distinct branches.  This is a general phenomenon
which applies to any homogeneous polynomial on $\Symn$ which is  G\aa rding hyperbolic
with respect to the identity. (See [HL$_{7,8}$] and \DD.3.4 below.)

\smallskip
\noindent
$\bullet$  \ {\bf (The Special Lagrangian Potential Equation).}  This equation $F(c)$, given in 
\BB.2.1(d), can be treated for all values of $c$ and has the nice feature that 
$\ft(c) = F(-c)$.

\smallskip\noindent
$\bullet$  \ {\bf (Geometrically Determined Subequations -- Calibrations).} These are subequations
determined by a compact subset $\GG$ of the Grassmann bundle of tangent $p$-planes
 by requiring that $\tr_W (\Hess u) \geq0$ for all $W\in \GG$. These include many interesting
 examples, including the subequations in calibrated geometry discussed at the outset.
 It also includes a new polynomial differential equation in Lagrangian geometry (see \BB.1.11(d)).
 Incidentally, this equation has branches whose study is a non-trivial application of the
 G\aa rding theory above.

\smallskip\noindent
$\bullet$  \ {\bf (Equations Involving the Principal Curvatures of the Graph and the $k$-Laplacian).}
For all such invariant equations on $G/H$, Theorem \HH.1.3 gives  existence (but not uniqueness).
 Strict boundary convexity is easily computable (see [HL$_6$, \S 17] for example).
Existence holds  on {\sl all} domains  for the $k$-Laplacian 
$|\nabla u|^2\D u + (k-2) (\n u)^t (\Hess\,u) (\n u) =0$,
when $1< k\leq \infty$  and   
when $k=1$ on mean-convex domains, where uniqueness fails catastrophically.

 \bigskip 
 A fundamental  point is that all such  equations can be carried over to any
 riemannian manifold with an appropriate (not necessarily
 integrable!) reduction of structure group.  This is done  by using the {\bf riemannian hessian} 
given in \S 8.2.     Theorem \HH.1.2 can then be applied, and we obtain the following corollary.
Let $\bbf$  and  $\bbm$ be constant coefficient subequations in $\rn$ with invariance group $G$.

\Theorem{\HH.2.2}  {\sl
Let $F$ be a subequation with monotonicity cone $M$ canonically determined
by   $\bbf$ and $\bbm$ on a riemannian manifold $X$ with a topological
$G$-structure. Let $\O\ss\ss X$ be a domain with smooth boundary which is both
$F$ and $\ft$ srictly convex. Assume there exists a strictly $M$-subharmonic function
on $\ob$.  Then the Dirichlet Problem for $F$-harmonic functions is uniquely solvable 
for all  $\vf\in C(\bo)$.
}

\medskip\noindent
$\bullet$  \ {\bf (Universal Riemannian Subequations).} Any constant coefficient
subequation $\bbf$ which in invariant under the natural action of O$(n)$
carries over directly to any riemannian manifold., and Theorem \HH.2.2 applies.  This includes most
of the examples above.

\smallskip\noindent
$\bullet$  \ {\bf (Universal Hermitian Subequations).} A constant coefficient
subequation $\bbf$  invariant under U$(n)$
carries over  to  any almost complex hermitian manifold.  There is a quaternionic
analogue. More generally, we have:

\smallskip\noindent
$\bullet$  \ {\bf (Equations on Manifolds with $G$-Structure).} A constant coefficient
subequation $\bbf$  invariant under a subgroup $G\ss {\rm O}(n)$ carries over to any
manifold equipped with a topological $G$-structure (see \HH.2.1).  This includes manifolds with
topological (or quasi) calibrations based on any fixed form in $\L^p\rn$.  Even the extreme
case $G=\{e\}$ is interesting here. An $\{e\}$-structure is a topological trivialization of $TX$.  
It transplants {\sl every} constant coefficient equation to $X$, and Theorem \HH.2.2 applies.
This holds, for example,
for  every orientable 3-manifold and every Lie group.
 
\bigskip

Theorem \HH.1.2 actually treats much more general equations on manifolds.
Affine  jet-equivalence gives  great flexibility to the result.

\medskip
\centerline
{\sl
Many variable-coefficient, inhomogeneous subequations on manifolds
}

\centerline
{\sl
can be transformed by local affine jet-equivalence
}

\centerline
{\sl
to  universally defined subequations, such as those in Theorem \HH.2.2,
}

\centerline
{\sl
while preserving the domains of strict boundary convexity.
}

\medskip

\smallskip\noindent
$\bullet$  \ {\bf (Calabi-Yau-Type Equations).} This is a good example of the power of affine jet equivalence. 
It applies to treat equations of type 
$
\left(  i \partial\overline{\partial} u  + \o \right)^n =  F(x,u) \o^n
$
on {\sl almost complex} hermitian manifolds, where $F>0$ is non-decreasing in $u$.
See \CC.2.8.

\smallskip\noindent
$\bullet$  \ {\bf (Inhomogeneous Equations).}  Many homogeneous equations can be transformed into
inhomogeneous equations by affine jet equivalence.  For example, from the
$k^{\rm th}$ branch of the Monge-Amp\` ere equation one can obtain:  $\l_k(\Hess u) = f(x)$
for any continuous function $f$.
See \CC.2.7.

\smallskip\noindent
$\bullet$  \ {\bf (Obstacle Problems).}  The methods here apply also to the Dirichlet Problem
with an Obstacle.  In this case not all boundary data are allowed. They are constrained by the
obstacle function.  This is another example of an inhomogeneous equation.  See \S \HH.6.

\smallskip\noindent
$\bullet$  \ {\bf (Parabolic Equations).} Each of these subequations has a parabolic cousin,
where existence and uniqueness results are generally stronger.  See \HH.5.

 \bigskip

 For any  subequation $F$ on a manifold $X$, one has the very natural
 \medskip
 \noindent
 {\bf Restriction Question:}
  {\sl When is the restriction of an $F$-subharmonic function on $X$
 to a submanifold $j:Y\ss X$,  a  $j^*(F)$-subharmonic function on $Y$?}
 
 \medskip
 
 For $C^2$-functions, this always holds,  and if fact defines the induced
subequation $j^*F$. 
However,  it is important and
  non-trivial for general upper semi-continuous subharmonics.
 There are several restriction results established in [HL$_9$].  They  are 
 relevant to calibrated and riemannian geometry.  Sometimes they lead
 to  characterizing $F$-subharmonics in terms of their restrictions to 
 special submanifolds. 
 
  An important case of this latter phenomenon occurs in almost complex
  manifolds. The ``standard'' way of defining plurisubharmonic functions is to require that the 
  restrictions to (pseudo) holomorphic curves are subharmonic. There also exists
   an intrinsic subequation,  whose subharmonics agree with
  the standard plurisubharmonic functions in the integrable case.
Via the restriction theorem, these two definitions have been shown to
 agree on any almost complex manifold [HL$_{10}$].
 
 There is also the notion of a plurisubharmonic distribution on a general
 almost complex manifold.  Nefton Pali [P]  has shown that  those
 which are  representable by continuous
 $[-\infty,\infty)$-valued functions are of the type above, and he conjectured
 that this should be true generally.  This leads to another topic.
 
 For convex subequations which are ``second-order complete'', 
 a Strong Bellman Principle can be applied.  It enables one to
 prove that distributionally $F$-subharmonic functions correspond in a 
 very precise sense to the  upper semi-continuous  $F$-subharmonic functions
 considered here. This is done in [HL$_{13}$].  Such arguments apply to prove
 the Pali Conjecture [HL$_{10}$].

\bigskip

\noindent
{\bf Some Historical Notes.}  There is of course a vast literature on
the {\sl principal} branches of $\cp$ and $\cp^\bbc$  of
the real and complex Monge-Amp\`ere equations.  Just to mention a few of the  historically
significant contributions beginning with Alexandrov:  [Al], [Po$_*$], [RT], [B], [W], 
[TU], [CNS$_*$], [CKNS], [BT$_*$], 
[HM],  [S$_1$],  [CY$_*$], and [Yau].
Quaternionic subharmonicity and the principal branch $\cp^\bbh$ of the quaternionic
Monge-Amp\`ere equation have been studied in [A$_*$] and [AV].
On compact complex manifolds without boundary, viscosity solutions  to equations of the form 
$
\left(  i \partial\overline{\partial} u  + \o \right)^n =  e^\vf v,
$
where $v>0$ is a given smooth volume form,  were studied in  [EGZ]. By establishing a 
comparison principle they obtain existence and uniqueness of solutions in important
borderline cases ($\o\geq0$, $v\geq0$ with $\int v>0$), and also show that these are  the
unique solutions  in the pluripotential sense.

The parabolic form of the 1-Laplacian gives rise to mean curvature flow by the 
level set method.  Some of the  interesting  results on this topic (see [ES$_*$], [CGG$_*$], [E], [Gi])
can be carried over from euclidean space to the riemannian setting by the 
methods of [HL$_6$].

The first basic work on the Dirichlet Problem for the convex branches of the
Special Lagrangian potential equation appeared in [CNS$_2$], and 
there are further results by Yuan [Y], [WY].

In  [AFS]  and [PZ] standard viscosity theory has been  adapted
to  riemannian manifolds by using the distance function, parallel translation,
Jacobi fields, etc.   For the problems considered  here  this machinery
in not necessary.  

In [S$_{2,3,4}$],    Z. Slodkowski  developed an axiomatic perspective on
generalized subharmonic functions, and addressed the Dirichlet Problem  in this context. 
He studied certain invariant ``pseudoconvex classes''  of functions  on euclidean space and
complex homogeneous spaces.
There is a version of duality which plays an important role in his  theory.
It is  formulated   differently from the one here.  However, in  the cases of overlap
 the two notions of duality are equivalent.
Interestingly, his results are used to prove a duality theorem 
for complex interpolation of normed spaces 
[S$_{5}$]

\bigskip
 \noindent
 {\bf Concerning Regularity.} In this paper there is no serious discussion of regularity for solutions 
 of the Dirichlet Problem. Indeed,  with the level of degeneracy allowed here, no regularity
 above continuity can be claimed generally.  Consider $u_{xx}=0$ in $\bbr^2$ for example.
(See also [Po$_1$] and [NTV] and references therein.) 
A good account of regularity results can be found in [E].
A general exposition  of viscosity methods and results appears in [CIL] and [C].

\bigskip
 \noindent
 {\bf Concerning $-\infty$.}  Our approach here is to steadfastly treat subsolutions
 from the point of view of classical potential theory. We  allow subsolutions
  ($F$-subharmonic functions) to assume the value $-\infty$,
 in contrast to standard viscosity theory where subsolutions are finite-valued.
  This has the  advantage
 of including basic functions, like the fundamental solution of the Laplacian, Riesz
 potentials, and $\log |f|$ with $f$ holomorphic, 
 into the class of subsolutions.  It  also allows the constant function
 $u\equiv-\infty$, which is crucial for the restriction theorems
 discussed in Chapter \II. This issue is not important for the Dirichlet Problem.

\vfill\eject


\noindent{\headfont \BB.\  Subequations -- a Geometric Approach.}
\medskip
 
The aim of this chapter is to present a geometric approach to  subequations, pioneered by Krylov [K].
This point of view clarifies and conceptually simplifies many aspects of the
theory. For transparency we begin with the basic case. 

\medskip
\noindent
{\bf \BB.1. Constant Coefficient Subequations in $\rn$.}
The 2-jets of functions on $\rn$ (i.e.,  Taylor polynomials of degree two) take
values in the vector space
$$
\bbj^2\ \equiv \ \bbr\times\rn\times\Symn \  \quad{\rm with\ traditional\ coordinates}\quad \ (r,p,A).
\eqno{(\BB.1.1)}
$$
\Def{\BB.1.1}  
A  {\sl second-order  constant coefficient subequation} on $\rn$ is a proper
closed subset $\bbf \ss \bbj^2$ satisfying the {\bf Positivity Condition}
$$
\bbf  + \cp \ \ss\ \bbf 
\eqno{(P)}
$$
and the {\bf Negativity Condition}
$$
\bbf  + \cn \ \ss\ \bbf 
\eqno{(N)}
$$
where 
$$
\cp \ \equiv\ \{(0,0,A)\in \bbj^2 : A\geq 0\}
\and
\cn \ \equiv \ \{(r,0,0) \in \bbj^2 : r\leq 0\},
$$
and the {\bf Topological Condition}
$$
\bbf \ =\ \overline{\Int \bbf}.
\eqno{(T)}
$$
We say $\bbf $ is {\sl pure second-order} if $\bbf =\bbr\times \rn\times \bbf _0$ for a 
closed subset $\bbf _0\ss\Symn$.  In this case only (P) is
required, since (N) is automatic and one can show that  (P)\ $\Rightarrow$\ (T).
Such subequations are often simply denoted  by the subset $\bbf_0$ of $\Symn$.

\Ex{\BB.1.2} Some basic pure second-order examples are:
\medskip
\noindent
(a) {\bf The Laplace Subequation:}  \smallskip

\centerline
{$\bbf _0 = \{A\in\Symn : \tr A\geq0\}$.}

\smallskip
\noindent
(b) {\bf The Homogeneous Monge-Amp\`ere Subequation:}
\smallskip
\centerline{\ $\bbf _0= \{A\in\Symn : A\geq0\}\ \cong\  \cp$.}

\smallskip
\noindent
(c) {\bf The $k^{\rm th}$ Elementary Symmetric Function Subequation:} 
\smallskip
\centerline{
$\bbf _0= \{A\in\Symn : \s_\ell(A)\geq0, 1\leq\ell\leq k\}$.}

 \smallskip

\noindent
(d) {\bf The Special Lagrangian Potential Subequation:} 
\smallskip
\centerline{
$\bbf _0= \{A\in\Symn :\tr ( \arctan\,A)\geq c\}$.}

 \smallskip

\noindent
(e) {\bf The Calabi-Yau Subequation:} (This is not pure second-order, but it is gradient-independent.)

\centerline{
$\bbf  = \{(r,p,A)\in\Symn : \tr(A+ I)\geq e^r \ {\rm and}\ A+I\geq0\}$.}

\Remark{\BB.1.3}
In $\bbc^n = (\bbr^{2n}, J)$ each of the examples above has a complex analogue
given by replacing $A$ with its hermitian symmetric part $A_\bbc \equiv \half(A-JAJ)$.
The same applies in quaternionic $n$-space $\bbh^n = (\bbr^{4n}, I,J,K)$ with 
$A$ replaced by   $A_\bbh \equiv {1\over 4}(A-IAI-JAJ-KAK)$.

\vfill\eject

\Def{\BB.1.4}  Given a constant coefficient subequation $\bbf $ on $\rn$, the {\bf dual}
subequation $\wt{\bbf}$ is defined by
$$
\wt{\bbf} \ \equiv\ \sim(- \Int \bbf ) \ =\ -(\sim \Int \bbf ).
$$

\Lemma {\BB.1.5} {\sl
\centerline{$\bbf $ is a subequation \ \ $\iff$\ \ $\wt{\bbf}$ is a subequation,\qquad\qquad\qquad\qquad\qquad\qquad}
\medskip
\noindent
and in this case
$$
\wt{\wt{\bbf}} \ =\ \bbf
\quad\and\quad
\wt{\bbf+J} \ =\ \wt \bbf - J
$$
for all $J\in \bbj^2$.}
\medskip
The proof can be found in [HL$_4$, \S 4].  In the examples above the dual subequations are easily 
computed in terms of the eigenvalues of $A$ (or $A_\bbc$, etc.). One finds that  the Laplace
subequation is self-dual ($\wt{\bbf}=\bbf $) but the others are generally not. Of particular interest is example (b)
where the dual of $\cp \equiv \{A\geq0\}$ is
$$
\cpt\ \cong \ \{A\in \Symn : {\rm at\ least\  one\ eigenvalue \ of \ } A\ {\rm is\ }\geq 0\}
\eqno{(\BB.1.2)}
$$

We now present  a concept of central importance which comes from viscosity theory [CIL].
For any manifold $X$, let  $\USC(X)$ denote  the set of 
upper semi-continuous functions $u:X\to  [-\infty, \infty)$. 
Given  $u\in \USC(X)$ and a point $x\in X$,  a {\bf test function for $u$ at $x$} 
is a $C^2$-function $\vf$ defined near $x$ so that 
$$
u\ \leq\ \vf
\and u(x) \ =\ \vf(x).
$$

\Def{\BB.1.6} Let $\bbf $ be a constant coefficient subequation on $\rn$ and fix an open set $X\ss\rn$.
A function $u\in\USC(X)$ is said to be {\bf $\bbf $-subharmonic} on $X$ if for each $x\in X$ and each
test function $\vf$ for $u$ at $x$, the  {\sl 2-jet}  (or {\sl total second derivative}) of $\vf$ satisfies
$$
J_x^2 \vf \equiv (\vf(x), (D\vf)_x, (D^2 \vf)_x) \ \in \ \bbf .
\eqno{(\BB.1.3)}
$$
It is important that this condition (\BB.1.3) is only required at points where test functions actually exist.
The set of such functions is denoted by $F(X)$.  

\medskip
It is striking that the space  $F(X)$ of $F$-subharmonics  shares  many of the 
 important properties enjoyed by classical subharmonic functions (see    \BB.3.1 below).
The $C^2$-functions $u\in F(X)$ are exactly those with $J^2_x u\in \bbf $ for all $x\in X$.
This basic fact requires the Positivity Condition (P) on $\bbf $.
Interestingly, the other properties in   \BB.3.1 do not require (P). 

In the example (a) we have the following (see [HL$_4$, Rmk. 4.9] and  [HL$_9$, Prop. 2.7]).

\Prop{\BB.1.7} {\sl
\smallskip

(i) \ \ $\cp(X)$ is  the set of convex functions on $X$.
\smallskip

(ii)\ \
$\cpt(X)$ is the set of subaffine functions on $X$.
}

\Def{\BB.1.8} A function $u\in \USC(X)$ is called {\bf subaffine} if for each compact subset $K\ss X$ 
and each affine function $a$, 
$$
u\ \leq\ a\quad{\rm on} \ \ \partial K\qquad \Rightarrow\qquad
u\ \leq\ a\quad{\rm on} \ \ K.
$$

Note that subaffine functions satisfy the maximum principle. 
In fact, for a pure second-order
subequations,  the subequation $\cpt$
is {\sl universal} for this property.
  That is,  if the functions in $\bbf(X)$ satisfy the maximum principle, then
$\bbf\ss \cpt$.
We note also that functions which are locally subaffine 
are globally subaffine, while the corresponding statement for functions 
satisfying the maximum principle is false.

\vfill\eject

\Def{\BB.1.9} Let $\bbf $ and $X$ be as in Definition \BB.1.6.
A function $u\in\USC(X)$ is said to be {\bf $\bbf $-harmonic} on $X$
if 
$$
u\ \in\ F(X)\and
-u\ \in\ \ft(X)
\eqno{(\BB.1.4)}
$$

Condition (\BB.1.4) implies that $u$ is continuous.  If $u$ is twice differentiable at a point
$x$, then (\BB.1.4) implies that 
$$
J_x^2 u \ \in \ \bbf \cap (-\wt \bbf) \ =\  \bbf \cap (\sim \Int \bbf) \ =\  \partial \bbf.
$$
 Thus if $\bbf $ is defined classically as the closure of
a set  $\{ f(r,p,A)> 0\}$ for a continuous function
$f:\bbj^2 \to \bbr$, then any  $u\in C^2(X)$ which is $\bbf $-harmonic satisfies the differential equation
$$
f(u, Du, D^2u)\ =\ 0 \qquad{\rm on\ \ } X,
$$
however, the converse is not always true.

\Note {\BB.1.10. (Branches)} It is instructive to consider the most basic of subequations, $\cp$. 
A $C^2$-function $u$ which is 
$\cp$-harmonic satisfies the homogeneous Monge-Amp\`ere equation 
$$
\det \left( D^2 u\right)\ =\ 0.
\eqno{(\BB.1.5)}
$$
However, $u$ is required to have the additional property of being convex (cf. Alexandroff  [Al]).
 (In the complex analogue $u$ is plurisubharmonic.) 

The equation (\BB.1.5) has other solutions corresponding to other ``branches''  of the
locus $\{\det A=0\}$, which can also be handled by this theory.  Given a symmetric matrix 
$A$, let $\l_1(A) \leq \l_2(A)\leq \cdots\leq \l_n(A)$ be the {\sl ordered} eigenvalues of $A$.
Since  $\det A= \l_1(A) \cdots \l_n(A)$, equation (\BB.1.5) can be split into branches
$$
\l_k\left(D^2 u\right) \ =\ 0.
\eqno{(\BB.1.5)_k}
$$
for $k=1,...,n$.
By monotonicity of eigenvalues, each   $\bL_k  \equiv \{\l_k\geq0\}$ is a subequation.
Interestingly, the dual of a branch is another branch:
$$
\wt{\bL_k} \ = \ \bL_{n-k+1}
$$

This   phenomenon of branches  occurs in many equations of geometric significance.

\Ex{\BB.1.11. (Geometrically Defined Subequations)} There is a large class of subequations which
arise naturally in our set-theoretic setting. 
Let $G(p,\rn)$ denote  the Grassmannian of $p$-planes in $\rn$. 
For each compact subset $\GG\ss G(p,\rn)$ we  define
the pure second-order  subequation
$$
\bbf(\GG) \ \equiv \ \{ A\in \Symn : \tr_W A\geq 0 \ \ {\rm for\ all\ } W\in \GG\}
\eqno{(\BB.1.6)}
$$
with dual
$$
\wt{ \bbf (\GG)} \ =\ \{ A\in \Symn : \tr_W A\geq 0 \ \ {\rm for\ some\ } W\in \GG\}
$$
The $\bbf(\GG)$-subharmonic functions are called $\GG$-{\sl plurisubharmonic}.
This terminology is justified by the following. Let $X\ss \rn$ be an open set.

\Theorem{\BB.1.12} {\sl 
A function $u\in \USC(X)$ is $\GG$-\psh if and only if for every affine $\GG$-plane $L$
the restriction $u\bigr|_{X\cap L}$ is subharmonic for the standard Laplacian on $L$.
 The same statement holds with the affine $\GG$-planes  expanded to include
 all  minimal $\GG$-submanifolds of $X$.
 (A $\GG$-submanifold is one whose tangent planes are elements of  $\GG$).
 }

\vfill\eject

This follows from a Restriction Theorem in [HL$_9$], which  is discussed in Chapter \II.
\medskip
 \noindent
 (a) \ \ $\GG = G(1,\rn)$:  In this case $\bbf(\GG) = \cp$ and the $\GG$-\psh functions are the
classical convex functions, i.e., those which are convex on affine lines.
 
 \medskip
 \noindent
 (b) \ \ $\GG = G_\bbc(1,\bbc^n) \ss G(2,\bbr^{2n})$ the set of complex lines in $\bbc^n$:  
 In this case $\bbf(\GG) = \cp^\bbc$ (see \DD.3.1), and the $\GG$-\psh functions are the
 standard plurisubharmonic functions, i.e., those which are subharmonic on complex lines.

\medskip
 \noindent
 (c) \ \ $\GG = G(p,\rn)$:  Here  the $\GG$-\psh functions are the
standard $p$-plurisubharmonic functions, i.e., 
 those which are subharmonic on
 affine $p$-planes.   This subequation 
 has the feature that 
 each $p$-\psh function is also $\GG$-\psh for every closed $\GG\ss G(p,\rn)$.
 The   analogue $\GG = G(p,\bbc^n)$ in the complex case
 plays a  role in analysis in several complex variables.

 \medskip

  The $\GG$-harmonic functions in these cases are
 viscosity solutions to differential equations which are O$(n)$  (or U$(n)$) invariant
  polynomials  in the variables  $D^2u$. Each of these equations has branches
 which will be discussed further  in \DD.3.1 and \DD.3.2 below.

\medskip
 \noindent
 (d) \ $\GG = {\rm LAG} \ss G(n, \bbr^{2n})$ the set of Lagrangian planes in $\bbc^n=\bbr^{2n}$:
   In this case  the ${\rm LAG}$-\psh functions are relatively new and interesting. The corresponding
   harmonics are viscosity solutions to a differential equation  which is a U$(n)$-invariant polynomial
    in the variables   $D^2u$ (see [HL$_{14}$]). This equation also has branches.

\medskip
Many important  examples come directly from the theory of calibrations. A {\sl  parallel calibration}
in $\rn$ is a  constant   coefficient $p$-form whose restriction satisfies $\pm\vf|_W \leq { vol}_W$
 for all oriented $p$-planes $W$. 
For such a $\vf$, we define $\GG\equiv G(\vf)$ to be the set of $W\in G(p,\rn)$ such that 
$|\vf|_W |=   { vol}_W$.
In this case $G(\vf)$-submanifolds (or simply $\vf$-submanifolds) are automatically minimal. 
When $\vf=\o$ is the K\"ahler form in 
$\bbc^n$, we recover case (b) above, where the $\o$-submanifolds are the holomorphic curves.
(This carries  over to any symplectic 
manifold $(X,\o)$ with a compatible almost complex structure in the sense of Gromov [Gr].)
The $G(\vf)$-\psh (or simply $\vf$-plurisubharmonic) 
functions are  essentially  {\bf dual} to the $\vf$-submanifolds (see [HL$_{2,3}$]),
and they provide calibrated geometry with new tools from conventional analysis.

\medskip
 \noindent
 (e) \  $\GG = G(\vf) = {\rm SLAG}  \ss G(n, \bbr^{2n})$ where $\vf= {\rm Re}(dz_1\wedge \cdots\wedge dz_n)$ is the Special Lagrangian Calibration (cf. [HL$_{1}$]). The notions of Special Lagrangian submanifolds and of
 SLAG-\psh  and SLAG-harmonic functions carry over to any Ricci-flat K\'ahler manifold (cf. [HL$_{1}$]). 
 The SLAG-subvarieties play a central role in the conjectured differential-geometric interpretation of mirror symmetry presented in  [SYZ$_{1,2}$].

\medskip
 \noindent
 (f) \  $\GG = G(\vf)  \ss G(3, \bbr^{7})$ where $\bbr^7={\rm Im}{\bf O}$ is   the imaginary octonions
 and $\vf(x,y,z) \equiv \bra{x\cdot y} z$ is the {\bf associative} calibration. There is a rich geometry of associative submanifolds, and an abundance of $\vf$-\psh and $\vf$-harmonic functions. 
 The same applies to the {\bf coassociative} 
 calibration $\psi= * \vf$.  Both calibrations make sense on any 7-manifold with G$_2$-holonomy.

\medskip
 \noindent
 (g) \  $\GG = G(\Phi)  \ss G(4, \bbr^{8})$ where $\bbr^8={\bf O}$, the  octonions, 
 and $\Phi(x,y,z,w) \equiv \bra{x\times y\times z} w$ is the {\bf Cayley} calibration. 
 There is a rich geometry of Cayley  submanifolds, and an abundance of $\Phi$-\psh 
 and $\Phi$-harmonic functions.  All this carries over to  any 8-manifold with Spin$_7$-holonomy.

\medskip
\noindent
{\bf Note.}  While the $\vf$-harmonic functions in examples (e), (f) and (g) are of basic interest
in calibrated geometry, they appear {\bf not} to satisfy any polynomial equation in $u, Du$ and $D^2u$.
This is one justification for the the  approach to subequations adopted here.


\vskip.3in
\noindent
{\bf \BB.2. Subequations on General  Manifolds.}
Suppose now that $X$ is a smooth manifold of dimension $n$.
The natural setting for second-order differential equations on $X$
is the bundle of {\bf 2-jets} of functions on $X$.  This is the bundle $J^2(X) \to X$ whose
fibre at $x\in X$ is the quotient $J_x^2(X) = C^\infty_x/ C^\infty_{x,3}$ of germs 
of smooth functions at $x$ modulo those which vanish to order 3 at $x$.

Restriction from 2-jets to 1-jets gives  a basic short exact sequence
$$
0\ \arr\ \Sym(T^*X) \ \arr\ J^2(X)\ \arr\ J^1(X) \ \arr\ 0
\eqno{(\BB.2.1)}
$$
where $\Sym(T_x^*X)$ embeds into $J^2_x(X)$ as the 2-jets of functions having
a critical value zero at $x$. The dual exact sequence is
$$
0\ \arr\  J_1(X) \ \arr\ J_2(X)\ \harr {\s}{}\  \Sym(TX)  \ \arr\ 0.
\eqno{(\BB.2.2)}
$$
Sections of $J_k(X)$ are linear differential operators of degree $\leq k$ on $X$,
and $\s$ is the {\sl principal symbol map} on operators of  degree $2$.
 
There are two important, intrinsically defined subbundles of $J^2(X)$ which correspond to the 
subspaces $\cp$ and $\cn$ in Definition \BB.1.1 , namely:
$$
\cp\ \equiv\ \{A\in \Sym(T^*X) : A\geq0\}
\and
\cn \ \equiv \ \{{\rm 2-jets \ of\  \ constant\ functions\ } \leq0\}.
$$

\Def{\BB.2.1} A {\sl  subequation} of order $\leq2$ on $X$ is a closed subset
$F\ss J^2(X)$  satisfying  (under fibre-wise sum) the {\sl Positivity Condition}:
$$
F + \cp \ \ss\ F,
\eqno{(P)}
$$
the {\sl Negativity Condition}:
$$
F + \cn \ \ss\ F,
\eqno{(N)}
$$
and the {\sl Topological Condition}:
$$
 (i)\ \ F\ =\ \overline{\Int F}, \qquad (ii)\ \ F_x\ =\ \overline{\Int F_x}, \qquad 
 (iii)\ \ \Int F_x\ =\     (\Int F)_x
 \eqno{(T)} 
 $$
\smallskip
 \noindent
 where $\Int F_x$ denotes interior with respect to the fibre.

\medskip

Note that $\cp$ is {\sl not} a subequation. However, 
when discussing  pure second-order subequations,  it is sometimes 
used as an abbreviation for  $\bbr\times\rn\times\cp$, which is a subequation.
(see \BB.1.1 and \BB.1.2).

\Remark {\BB.2.2.  (Splitting the 2-Jet Bundle)} Let $\n$ be a torsion-free connection
on $X$. Then each $u \in C^2(X)$   has an associated hessian $\Hess \,u \in\G ( \Sym(T^*X))$
defined on vector fields $V,W$ by
$$
(\Hess\, u)(V,W) \ =\ V W u -  W V u - (\n_V W) u.
\eqno{(\BB.2.3)}
$$
Since $\n_V W -\n_W V = [V,W]$, one easily sees that $\Hess\, u$ is a symmetric tensor.
If $X$ is riemannian and $\n$ is the Levi-Civita connection, then $\Hess\, u$ is called
the {\sl riemannian hessian} of $u$.  

The hessian in (\BB.2.3) depends only on the 2-jet of $u$ at each point, and so it gives
a splitting of the short exact sequence (\BB.2.1).  That is, we can write
$$
J^2(X) \ =\ \bbr\oplus T^*X \oplus \Sym(T^*X)
\eqno{(\BB.2.4)}
$$
by the association
$$
 J^2_x u = (u(x), (du)_x, \Hess_x u).
$$

\Remark {\BB.2.3.  (Universal Subequations)} Each of the subequations given in 
Example \BB.1.2 carries over to any riemannian manifold $X$ by using the splitting
(\BB.2.4) (determined by the riemannian hessian). For instance, Example 2.1.2(a) gives the
Laplace-Beltrami operator. More generally, any 
constant coefficient subequation $\bbf \ss \bbj^2$ which is invariant under the
action of the group O$(n)$,  transplants to every riemannian manifold.
In the case of $\bbc^n=(\bbr^{2n}, J)$, each U$(n)$-invariant subequation
transplants to every hermitian almost complex manifold.

There is, in fact, a very general principle:
\medskip
\centerline
{
\sl Let $\bbf \ss \bbj^2$  be a constant coefficient subequation which is invariant } 
\centerline
{
\sl under a subgroup $G\ss $ O$(n)$ acting naturally on $\bbj^2$.
} 
\centerline
{
\sl Then $\bbf$ carries over to a subequation $F$ on every manifold $X$ with  a topological $G$-structure.} 

\medskip
\noindent

The reader is referred to [HL$_6$] and \S \HH.2 below for definitions and many examples.

The concepts of the previous section now carry over to this general  setting.

\bigskip

\Def{\BB.2.4}  Given a  subequation $F\ss J^2(X)$, the {\bf dual}
subequation $\ft$ is defined by
$$
\ft \ \equiv\ \sim(- \Int F) \ =\ -(\sim \Int F).
$$

\Lemma {\BB.2.5} {\sl
\centerline{$F$ is a subequation \ \ $\iff$\ \ $\ft$ is a subequation,\qquad\qquad\qquad\qquad\qquad\qquad}
\medskip
\noindent
and in this case
$$
\wt{\ft} \ =\ F
\and 
\wt{F+S}\ =\ \ft-S
$$
for any section $S$ of $J^2(X)$.
}
\medskip
The proof can be found in [HL$_6$ \S 3].   The dual  of a universal subequation 
associated to  $\bbf\ss\bbj^2$ is the universal subequation associated to  $\wt\bbf$.
As before we have the following.

\Def{\BB.2.6} Let $F$ be a   subequation on a manifold $X$.
A function $u\in\USC(X)$ is said to be {\bf $F$-subharmonic} on $X$ if for each $x\in X$ and each
test function $\vf$ for $u$ at $x$, 
$$
J_x^2 \vf \equiv (\vf(x), (D\vf)_x, (D^2 \vf)_x) \ \in \ F.
\eqno{(\BB.2.5)}
$$
The set of such functions is denoted by $F(X)$.

\Def{\BB.2.7}  Let $F$ be a   subequation on a manifold $X$.
A function $u\in\USC(X)$ is said to be {\bf $F$-harmonic} on $X$
if 
$$
u\ \in\ F(X)\and
-u\ \in\ \ft(X)
\eqno{(\BB.2.6)}
$$

As before, positivity ensures that a function $u\in C^2(X)$ is $F$-subharmonic on $X$ iff $J^2_x u\in F$  for all $x$, 
and it is $F$-harmonic iff
$$
J_x^2 u \ \in \ \partial F \fa x.
$$


\vskip .3in
\noindent
{\bf \BB.3. Properties of $F$-Subharmonic Functions.}
The  $F$-subharmonic  functions  share many of the important properties
of classical subharmonic functions.  

\Theorem{\BB.3.1.  (Elementary  Properties of  F-Subharmonic Functions)}
{\sl
Let $F$ be an arbitrary  closed subset of $J^2(X)$.
\medskip

\item{(i)}  (Maximum Property)  If $u,v \in F(X)$, then $w=\max\{u,v\}\in F(X)$.

\medskip

\item{(ii)}     (Coherence Property) If $u \in F(X)$ is twice differentiable at $x\in X$, then $J_x^2 u\in F_x$.

\medskip

\item{(iii)}  (Decreasing Sequence Property)  If $\{ u_j \}$ is a 
decreasing ($u_j\geq u_{j+1}$) sequence of \ \ functions with all $u_j \in F(X)$,
then the limit $u=\lim_{j\to\infty}u_j \in F(X)$.

\medskip

\item{(iv)}  (Uniform Limit Property) Suppose  $\{ u_j \} \ss F(X)$ is a 
sequence which converges to $u$  uniformly on compact subsets to $X$, then $u \in F(X)$.

\medskip

\item{(v)}  (Families Locally Bounded Above)  Suppose $\cf\subset F(X)$ is a family of 
functions which are locally uniformly bounded above.  Then the upper semicontinuous
regularization $v^*$ of the upper envelope 
$$
v(x)\ =\ \sup_{f\in \cf} f(x)
$$
belongs to $F(X)$.

}
\medskip

A proof can be found,  for example,  in Appendix B in [HL$_6$]. 
For parts (i) and (ii),
even the closure hypothesis on $F$ can be weakened (op. cit.).

\vfill\eject


\noindent{\headfont \CC.\  Jet Equivalence of Subequations.}
\medskip
Many important nonlinear equations that occur in geometry can be 
transformed locally to constant coefficient equations.  This technique allows one to 
apply standard arguments from viscosity theory to prove local comparison results.

\medskip
{\bf \CC.1.  Affine Automorphisms of the Jet Bundle $J^2(X)$.}
 The transformations we shall use are the affine  automorphisms of $J^2(X)$ which we now 
 introduce.  To begin, note that there is a canonical direct sum decomposition
 $$
 J^2(X) = \bbr\oplus J^2_{\rm red}(X)
 \eqno{(\CC.1.1)}
 $$
 where the trivial $\bbr$-factor corresponds to the value of the function.  For the reduced 2-jet 
 bundle there is a short exact sequence
 $$
  0\ \arr\ \Sym(T^*X)\ \arr\ J^2_{\rm red}(X)\ \arr\ T^*X\ \arr\ 0
  \eqno{(\CC.1.2)}
 $$
coming from (\BB.2.1) above.

\Def{\CC.1.1} A linear isomorphism of $J^2(X)$ is an  {\bf automorphism} 
if, with respect to the splitting (\CC.1.1) it has the form ${\rm Id}\oplus \Phi$
where $\Phi:  J^2_{\rm red}(X)\to  J^2_{\rm red}(X)$ has the following properties.
We first require that
$$
\Phi(\Sym(T^*X)) = \Sym(T^*X), 
\eqno{(\CC.1.3)}
$$
 so by (\CC.1.2) there is an induced bundle automorphism
$$
g=g_{\Phi}:T^*X \ \arr\ T^*X.
\eqno{(\CC.1.4)}
$$
We further require that there exist a second bundle automorphism 
$$
h=h_{\Phi}:T^*X \ \arr\ T^*X
\eqno{(\CC.1.5)}
$$
such that on $\Sym(T^*X)$, $\Phi$ has the form $\Phi(A) = hAh^t$, i.e., 
$$
\Phi(A)(v,w) \ =\ A(h^tv,h^tw)    \qquad {\rm for\ } v,w\in TX.
\eqno{(\CC.1.6)}
$$

  The automorphisms of $J^2(X)$ form a group. They are the sections 
of the bundle of groups Aut$(J^2(X))$ whose fibre at $x\in X$ is the group of automorphisms
of $J_x^2(X)$ defined by (\CC.1.3) - (\CC.1.6) above. See [HL$_6$, \S 6.2] for this and
the following.

\Prop{\CC.1.2} {\sl With respect  to any  splitting 
$$
J^2(X) \ =\ \bbr\oplus T^*X \oplus \Sym(T^*X)
$$
of the short exact sequence (\BB.2.1), a bundle automorphism has the form
$$
\Phi(r, p,  A) \ =\ (r, gp, hAh^t + L(p))
\eqno{(\CC.1.7)}
$$
where  $g, h:T^*X\to T^*X$ are bundle isomorphisms
and $L$ is a smooth section of the bundle 
$\Hom (T^*X, \Sym(T^*X))$. }

\Ex{\CC.1.3}  Given a local coordinate system $(\x_1,...,\x_n)$ on an open set
$U\ss X$, the  {\sl canonical trivialization}
$$
J^2(U)\ =\ U\times\bbr\times\rn\times\Symn
\eqno{(\CC.1.8)}
$$
is determined by $J^2_xu  = (u,Du,D^2u)$ 
where $Du = (u_{\x_1}, ... , u_{\x_n})$ and $D^2u = (\!( u_{\x_i \x_j} )\!)$
evaluated at the point $\x(x)\in\rn$.
With respect to this splitting, every automorphism is of the form
$$
\Phi(u, Du, D^2u) \ =\ (u, \ g Du, \ h\cdot D^2u\cdot h^t + L(Du))
\eqno{(\CC.1.9)}
$$
where $g_x, h_x \in \GL_n$ and $L_x : \rn\to \Symn$ is linear for each point $x\in U$.

\Ex{\CC.1.4} The trivial 2-jet bundle on $\rn$ has fibre
$$
 \bbj^2  = \bbr\times\rn\times \Symn.
 $$
with automorphism group
 $$
{\rm Aut}(\bbj^2) \ \equiv\  {\rm GL}_n\times {\rm GL}_n\times \Hom(\rn, \Symn)
 $$
 where the  action is  given by
 $$
 \Phi_{(g,h,L)}(r,p,A) \ =\ (r, \ gp, \ hAh^t +L(p)).
 $$
Note that  the group law is
 $$
(\bar g, \bar h, \bar L)\cdot  (g, h,L)\ =\ (\bar g g,\ \bar h h, \ \bar h L \bar h^t + \bar L\circ g)
 $$

 Automorphisms at a point, with $g=h$,  appear naturally when one considers the action of diffeomorphisms.
 Namely, if $\vf$ is a diffeomorphism
fixing a point $x_0$, then in local coordinates (as in Example \CC.1.3 above) the right action on $J_{x_0}^2$, induced
by the pull-back $\vf^*$ on 2-jets,  is an automorphism.

\Remark{\CC.1.5}  Despite this last remark, automorphisms of the 2-jet bundle $J^2(X)$,
even those with $g=h$,  have little
to do with global diffeomorphisms or global changes of coordinates. In fact an automorphism radically restructures
$J^2(X)$ in that the image of an integrable section (one obtained by taking $J^2u$ for a fixed smooth function $u$ on $X$) is essentially never integrable.

\medskip

 The automorphism group ${\rm Aut}(J^2(X))$ can be naturally extended
by the fibre-wise translations.  Recall that the group
of affine transformations of a vector space $V$ is the product ${\rm Aff}(V) = \GL(V) \times V$
acting on $V$ by $(g,v)(u) = g(u)+v$.  The  group law is $(g,v)\cdot (h,w) = (gh, v+g(w))$.
There is a short exact sequence 
\smallskip
\centerline
{
$
0\to V\to {\rm Aff}(V) \harr{\pi} {} \GL(V) \to \{I\}.
$
}

\Def{\CC.1.6}  
The {\bf affine automorphism group} of $J^2(X)$ is   the space of smooth sections of
$$
\pi^{-1}\{  {\rm Aut}(J^2(X)) \}) \ \ss\ {\rm Aff}(J^2(X))
$$
where $\pi$ is the surjective bundle map
$\pi:{\rm Aff}(J^2(X)) \to \GL(J^2(X))$.  \medskip

Note that any affine automorphism  can be written in the form
$$
{ \Psi}\ =\ \Phi   + S
\eqno{(\CC.1.10)}
$$
where $\Phi$ is a (linear) automorphism and $S$ is a section of the bundle $J^2(X)$.


\vskip .3in
\noindent
{\bf \CC.2.  Jet-Equivalence.}

\Def{\CC.2.1}  Two subequations $F, F'\ss J^2(X)$ are said to be   {\bf jet-equivalent}  if 
there exists an automorphism  $\Phi:J^2 (X)\to  J^2(X)$ with 
$\Phi(F)=F'$.  If this holds for an affine automorphism $\Psi = \Phi +S$, 
they are said to be {\bf affinely jet-equivalent}.

\Remark{ \CC.2.2}  A  jet-equivalence  $\Phi: F\to F'$ does not take $F$-subharmonic functions
to $F'$-subharmonic functions.   In fact as mentioned above, for  $u \in C^2$,  
$\Phi(J^2 u)$ is almost never the 2-jet of a 
function. It happens if and only if  $\Phi(J^2 u) = J^2u$.
Nevertheless,  if $\Psi = \Phi+S$ is an affine  automorphism of $J^2(X)$ and $F\ss J^2(X)$ is a closed set,
then
\smallskip
\centerline
{
$F$  is a subequation \qquad$\iff$ \qquad
$\Psi(F)$ is   a subequation,
}
\smallskip 
\noindent
and furthermore, by \BB.2.5, 
$$
\wt {\Psi(F)} \ =\ \Phi(\ft) - S,
$$
which is  basic in establishing comparison.

\Def{\CC.2.3}  We say that a  subequation $F\ss J^2(X)$  {\sl is locally  affinely
 jet-equivalent to a constant coefficient
subequation}  $\bbf$  if each point $x$ has a local coordinate neighborhood  $U$ such that,
in the canonical trivialization (\CC.1.8) of $J^2(U)$ determined by  those coordinates, $F$ is affinely jet-equivalent to the  constant coefficient subequation $U\times \bbf$.
 
 \medskip
 
 This concept is robust as shown by  the following lemma, whose proof is a straightforward calculation.
 
\Lemma{\CC.2.4} {\sl If $F$ is affinely jet-equivalent to $\bbf$ in some local coordinate
 trivialization of $J^2(U)$, then this is true in every local coordinate trivialization of $J^2(U)$.}
 \medskip
  
A basic reason for introducing this concept is the following (see [HL$_6$, Prop. 6.9]).
Let $X$ be a riemannian manifold with topological $G$-structure
for a subgroup $G\ss {\rm O}(n)$ (see (\HH.2.1)). 

\Prop {\CC.2.5} {\sl
Suppose that $F \ss J^2(X)$ 
is the subequation determined by a $G$-invariant constant coefficient subequation
$\bbf\ss \bbj^2$ (cf. 2.2.3 and 8.2).
  Then $F$ is locally jet-equivalent to $\bbf$ on $X$.
}

\Ex {\CC.2.6. (Universal Equations)} Basic examples come from universal riemannian equations ($G={\rm O}(n)$)
such as those given in Example \BB.1.2 (a), (b), (c),  and their complex analogues on almost
complex hermitian manifolds ($G={\rm U}(n)$) or the analogues on almost quaternionic
hermitian manifolds  ($G={\rm Sp}(n)$).  There are also the other branches of these
equations as discussed in Note \BB.1.10.  There are also the many geometric examples
coming from Lagrangian geometry and calibrated geometry which are discussed below.
\medskip

\Ex {\CC.2.7. (Inhomogeneous Equations)} 
Another important fact about  affine jet equivalence is that it can transform 
inhomogeneous equations into constant coefficient ones and vice versa.  
We present several illustrative examples here (and more in \HH.5).  They each have the structure
$F = \Psi(H)$, $H=\Psi^{-1}(F)$ where $F$ is a
pure second-order, universal riemannian subequation, and
$$
\Psi(A) \ \equiv\ hAh^t + S \ =\ \eta^2 A +S
$$
where $h(x) =\eta(x) {\rm Id}$, for $\eta:X\to \bbr$,  and $S:X \to \Sym T^*(X)$ is a translation term. 
\medskip
\noindent
(i)  Let  $F$ correspond to the $k^{\rm th}$ branch $\{\l_k(\Hess\,u)=0\}$ of the 
homogeneous Monge-Amp\`ere equation (see  \BB.1.10).
Taking $\eta\equiv 1$ and $S=-f(x) {\rm Id}$ shows that $F$ is  affinely jet-equivalent to the 
inhomogeneous equation
$$
\l_k(\Hess\, u)\ = \ f(x)
$$
for any smooth function $f$.  
This includes the Monge-Amp\`ere equation from 2.1.2(b) when written as 
$\l_{\rm min}(\Hess\,u)=0$. 

\medskip
\noindent
(ii) Let  $F$ correspond to the universal equation $\det(\Hess\,u) = 1$ with
$\Hess\,u \geq0$.  One can transform this to the inhomogeneous equation
$$
\det (\Hess\,u) = f(x) \quad {\rm with} \quad\Hess\, u\geq0
$$ 
for any smooth $f>0$ by choosing  $\eta = f^{-{1\over 2n}}$ and $S=0$.

\medskip
\noindent
(iii)  More generally, one can    transform the universal  subequation: 
 $\s_k( \Hess\,u)= 1$ and 
$\s_\ell(\Hess\,u)\geq 0$,  \ $1\leq \ell<k$, into the inhomogeneous  equation
$$
\s_k(\Hess\,u)\ =\ f(x)
\and
\s_\ell(\Hess\,u)\geq 0, \ \ 1\leq \ell<k
$$
for any smooth $f>0$ by  choosing $\eta = f^{-{1\over 2k}}$ and $S=0$.

\Ex{\CC.2.8.  (The Calabi-Yau Equation)}  Let $X$ be an almost complex hermitian manifold 
(a Riemannian U$_n$-manifold), and consider the subequation $F\ss J^2(X)$ determined by 
the euclidean subequation:  
$$
\det_\bbc\{A_\bbc +I\}\ \geq\ 1  \and  A_\bbc +I \ \geq \ 0
$$
where $A_\bbc \equiv \half(A-JAJ)$ is the hermitian symmetric part of $A$.
Let $f>0$  be a smooth positive function on $X$ and write $f=h^{-2n}$.
Consider the global affine automorphism of $J^2(X)$ given by
$$
\Psi(r,p,A)
\ =\ (r,p, h^2 A+ (h^2-1)I)
$$
and set $F_f =   \Psi^{-1}(F)$.  Then
$$
\eqalign
{
(r,p, A) \in F_f \quad &\iff \quad 
  \det_\bbc \{h^2(A_\bbc+I)\} \geq 1    \ \ {\rm and}\ \  h^2(A_\bbc+I) \geq 0 \cr
&\iff \quad 
  \det_\bbc \{(A_\bbc+I)\} \geq f    \ \ {\rm and}\ \ (A_\bbc+I) \geq 0  \cr
}
$$
so we see that the $F_f$-harmonic functions are functions $u$ with $\det_\bbc\{\Hess_\bbc u   + I \} = f$
 and $\Hess_\bbc u   + I \geq 0$ (quasi-plurisubharmonic).
If $X$ is actually a complex manifold of dimension $n$ with K\"ahler form $\o$, this last equation
can be written in the more familiar form
$$
\left( i \partial\overline{\partial} u  + \o \right)^n\ =\ f \o^n 
$$
\noindent
with $u$  quasi-plurisubharmonic.   

One can similarly treat the equation
$$
\left( i \partial\overline{\partial} u  + \o \right)^n\ =\ e^u f \o^n.
$$
or the same equation with $e^u$ replaced by any non-decreasing positive  function $F(u)$.

 \medskip
 
 The concept of affine jet equivalence plays a critical role in the study of intrinsically
 subharmonic functions on almost complex manifolds [HL$_{10}$].

\vfill\eject


\noindent{\headfont \DD.\  Monotonicity.}
\medskip
 
 A concept of fundamental importance here is that of a {\sl monotonicity cone} for a given subequation.
It is the key to establishing comparison and removable singularity theorems for equations which
are highly non-convex.

\vskip .3in
\noindent
{\bf \DD.1.  The Constant Coefficient Case.} Let $\bbf, \bbm\ss\bbj^2$   be  constant coefficient subequations.

\Def{\DD.1.1}  We say that $\bbm$ is a {\bf monotonicity subequation} for $\bbf$ if 
$$
\bbf+\bbm\ \ss\ \bbf.
\eqno{(\DD.1.1)}
$$
It follows directly from \BB.1.6 that the 
sum of an $\bbf$-subharmonic function and an $\bbm$-subharmonic function is
again $\bbf$-subharmonic, provided that one of them is smooth. 
Thus, the reader can see that monotonicity is related to approximation 
whenever $\bbm$ has the {\sl cone property}
$$
t\bbm\ \ss\ \bbm \qquad {\rm for}   \ \  0\leq t\leq 1.
$$
When this holds $M$ can be expanded so that each fibre is a convex cone 
with vertex at the origin (cf. \DD.1.4). Under  this added assumption 
$\bbm$  is called  a {\bf monotonicity cone}.

\Lemma{\DD.1.2} {\sl If $\bbm$ is a monotonicity cone for $\bbf$, then}
$$
\qquad\quad
\wt\bbf+\bbm\ \ss\ \wt\bbf \qquad{\sl and}
\eqno{(\DD.1.2)}
$$
$$
\bbf+\wt \bbf\ \ss\ \wt \bbm.
\eqno{(\DD.1.3)}
$$
\noindent
These elementary   facts are basic.  The first states that:
\medskip
\centerline
{  \sl
$\bbm$ is a monotonicity
cone for $\bbf$$\qquad\iff\qquad \bbm$ a monotonicity cone for $\wt\bbf$.
}
\medskip
\noindent
The second is the algebraic precursor to proving that:
\medskip
\centerline
{\sl
The sum of an $\bbf$-subharmonic function and an $\wt \bbf$-subharmonic
function 
}
\centerline{\sl is $\wt \bbm$-subharmonic. 
}
\medskip
 \noindent
If one of the two functions is smooth, this last result follows easily from the
definitions.   It is important,  because in most  cases, the $\wt\bbm$-subharmonic functions satisfy the following:

\medskip
\noindent
{\bf Zero Maximum Principle}: 
{\sl For any  compact set $K$ in the domain of $u$,}
$$
u\ \leq \ 0 \quad{\sl on} \ \ \partial K\qquad\Rightarrow\qquad
u\ \leq \ 0 \quad{\sl on} \ \ K.
\eqno{(ZMP)}
$$

 \Ex{\DD.1.3} The (ZMP) holds for $\wt\bbm$-subharmonic functions when
\medskip
\centerline
{
$
\bbm\ =\ 
\{(r,p,A)
\in\bbj^2 :  r\leq - \g |p|, \ \ p\in \cd\ \ {\rm and\ }\  A\geq0\}
$
}
\medskip
\noindent
where $\g>0$ and $\cd\ss\rn$ is a convex cone with non-empty interior (and vertex at 0).
See Appendix B for a proof and  further discussion of Examples.
Note incidentally that the smaller $M$ is, the easier it is to be a monotonicity cone for $F$, while
the larger $\wt M$ is, the harder it is  to satisfy (ZMP).

\Note{\DD.1.4} Associated to any   subequation $\bbf$ is the set
$\bbm_\bbf$   of all  $J\in \bbj^2$ such that   $ \bbf+tJ\ss \bbf$ for $0\leq t\leq 1$.
One checks easily that $\bbm_\bbf$ is a closed convex cone which satisfies (P) and (N).
Thus, if $\Int \bbm_\bbf \neq \emptyset$, it is the maximal monotonicity cone for $\bbf$.



\vfill\eject
\noindent
{\bf \DD.2.  The General  Case.} Let $F\ss J^2(X)$ be a subequation on a manifold $X$.

\Def{\DD.2.1}  A     {\bf monotonicity cone} for $F$  is a convex cone subequation $M\ss J^2(X)$
(each fibre  is a convex cone with vertex at the origin) satisfying the condition
$$
F+M\ \ss\ F
\eqno{(\DD.2.1)}
$$

\Lemma{\DD.2.2} {\sl If $M$ is a monotonicity cone for $F$, then}
$$
\qquad\quad
\wt F+ M\ \ss\ \wt F \qquad{\sl and}
\eqno{(\DD.2.2)}
$$
$$
F+\wt F\ \ss\ \wt  M.
\eqno{(\DD.2.3)}
$$

\Note{\DD.2.3}  Suppose $\bbf\ss\bbj^2$ is a constant coefficient subequation invariant under a
subgroup $G\ss {\rm O}(n)$.  Then $\bbm_\bbf$ is also $G$-invariant. Thus if $\Int \bbm_\bbf
\neq \emptyset$, it determines a monotonicity cone $M_F$  for every subequation $F$
canonically determined on any manifold with a topological $G$-structure (cf. Remark \BB.2.3).


\vskip .3in
\noindent
{\bf \DD.3. Examples. (Branches of Polynomial Equations)} Many subequations have naturally associated monotonicity cones.
The most basic case is the following.

\Ex {\DD.3.1. (Homogeneous Monge Amp\`ere Equations)} Let $K=\bbr, \bbc$ or $\bbh$
and let $K^n = \bbr^N$ for $N/n=1,2$, or 4.  Then any quadratic form $A\in\Sym(\bbr^N)$ has
a $K$-hermitian symmetric part $A_K$ defined in Remark \BB.1.3. 
Let $\l^K_1(A)\leq \cdots \leq \l^K_n(A)$ be the ordered eigenvalues of $A_K$
(where we ignore the natural multiplicities 2 in the complex case and 4 in the quaternion case).
Let 
$$
 \bL^K_k\ \equiv\ \{\l^K_k(A)\geq0\}
$$
denote the $k^{\rm th}$ branch of the homogeneous Monge-Amp\`ere equation (cf. Note \BB.1.10).
The dual subequation is $\wt\bL^K_k = \bL^K_{n-k+1}$. These subequations carry over to any 
riemannian manifold with orthogonal almost complex or quaternionic structures.  

The smallest, most basic branch is $\bL^K_1 = \{A^K\geq 0\} = \bbf(G(1, K^n))$, which will 
be denoted  by $ \cp^K$, $K=\bbr,\bbc$ or $\bbh$.
The monotonicity of ordered eigenvalues: 
$\l^K_k(A) \leq \l^K_k(A+P)$ for $P\in \cp^K$ implies that
$$
\bL^K_k + \cp^K\ \ss\ \bL^K_k,
$$
i.e., the top branch $\cp^K$ is a monotonicity cone for each branch $\bL^K_k$ of the   
Monge-Amp\`ere equation.

\Ex {\DD.3.2. ($p$-Convexity)} Fix $p$, 
$1\leq p\leq n$.  For each $A\in \Symn$ and each $p$-tuple $I = \{i_1< i_2< \cdots < i_p\}$, set $\l_I  (A)= 
\l_{i_1}(A) +\cdots +    \l_{i_p}(A)$.
Consider the second-order polynomial differential equation determined  by
$$
{\rm MA}_{p}(A) \ \equiv \ \prod_{I}\l_I(A)\ 
=\ \det \left \{   D_A : \L^p \bbr^n \to \L^p \bbr^n          \right\} \ =\ 0
$$
where $D_A$ denotes $A$ acting as a derivation on the exterior power $ \L^p \bbr^n$.
This equation splits into branches  
$
\bL_k(p)
$,
$k=1,..., {n\choose p}$, obtained by ordering the eigenvalues $\{\l_I(A)\}$.  The {\sl principle branch} $\bL_1(p)$,
which is denoted by
$$
\cp(p) \ \equiv\   \left\{A : \l_1(A)+\cdots+\l_p(A) \geq0  \right\} \ =\ \bbf(G(p,\rn)),
$$
is exactly the one  considered in \BB.1.11(c).  In particular, the  $\cp(p)$-subharmonic 
functions are just the     {\sl $p$-plurisubharmonic} functions -- those which are harmonic
on all affine p-planes. The monotonicity of eigenvalues 
shows that $\cp(p)$  is a monotonicity cone for
every branch of this equation, that is, 
$$
 \bL_k(p)  + \cp(p)  \ \ss\ \bL_k(p).   
 $$

More generally, let $K=\bbr, \bbc$ or $\bbh$ and, using the notation of \DD.3.1, set
$$
{\rm MA}_{p}^K(A) \ \equiv \ \prod_{I}\l_I^K(A).
$$
This defines a polynomial differential equation with principal branch
$\cp^K(p) = \bbf(G(p, K^n))$.  The  other  branches, obtained as above by ordering 
the   eigenvalues  $\{\l_I^K(A)\}$, are subequations
for which $\cp^K(p)$ is a monotonicity cone.

The cone $\cp(p)$ can be defined for any {\sl real} number $p$, $1\leq p\leq n$ by
$$
\cp(p) \ \equiv\   \left\{A : \l_1(A)+\cdots+\l_{[p]}(A) + (p-[p]) \l_{p+1}(A)\geq0  \right\}.
\eqno{(\DD.3.1)}
$$
This extension plays an important role in removable singularity theorems (see Section \FF.2  below).
We note that this extended $\cp(p)$ is the principal branch of the polynomial operator
${\rm MA}_p(A) = \prod(\l_I(A) + (p-[p])\l_k(A))$ where the product is over $|I|=[p]-1$ and $k\notin I$.

\Ex {\DD.3.3. ($\d$-Uniform Ellipticity)} A basic family of monotonicity subequations is given by
$$
\cp(\d) \ \equiv\ \{A\in\Symn : A\geq -\d \tr A \cdot I\}
$$
for $\d >  0$.  Any subequation $\bbf$, for which $\cp(\d)$ is a monotonicity cone, 
is  uniformly elliptic in the usual sense. This subequation is   the principal branch of the pure 
second-order polynomial differential equation:
$$
\prod_{i=1}^n (\l_k(\Hess\,u)+ \d \D u)\ =\ 0.
$$
This equation has $n$ branches 
$$
\l_k(\Hess\,u)+ \d \D u\ \geq\ 0 \qquad {\rm for \ } k=1,...,n, 
$$
and $\cp(\d)$ is a monotonicity cone for each of these branches, so in particular, 
each branch is uniformly elliptic.

This is easily generalized as follows.  Suppose $\bbf\ss\Symn$ is any pure second-order subequation.
Then for each $\d>0$, the {\sl $\d$-elliptic regularization $\bbf(\d)$} is defined by requiring that
$A+\d (\tr A)\cdot I \in\bbf$.  Now if $\bbm$ is a monotonicity cone for $\bbf$, it follows
immediately from the definitions that $\bbm(\d)$ is a monotonicity cone for $\bbf(\d)$.
Also, $\cp\ss \bbm$ implies that $\cp(\d)\ss \bbm(\d)$, which ensures that each $\bbf(\d)$ is
uniformly elliptic.

\def\PP{Q}

\Ex {\DD.3.4. (G\aa rding Hyperbolic Polynomials)}
The examples above, and several  below, fall into a general class of equations
where monotonicity cones appear naturally.    A   homogeneous
polynomial $\PP:\Symn\to\bbr$ of degree  $m$ is said to be {\sl G\aa rding hyperbolic
with respect to the identity} if $\PP(I)=1$ and for each $A\in\Symn$ the polynomial 
$q_A(t) \equiv \PP(tI + A)$ has $m$ real roots. Thus we can write
$$
\PP(tI+A)\ =\ \prod_{k=1}^m (t+\l_k(A))
$$
where the $\l_1(A)\leq  \cdots\leq \l_m(A)$ are  the
{\sl ordered eigenvalues} (the negatives of the roots)  of  $q_A(t)$.
Such a polynomial has $m$ branches
$$
{ \bL}_{\PP,k} \ \equiv\ \{ \l_k(A) \geq0\},\qquad k=1,...,m,
$$
which correspond to $m$ constant coefficient pure second-order subequations
in $\rn$.  The principal branch
$$
\bbm_\PP \ \equiv \ {\bL}_{\PP,1}
$$
is called the {\sl G\aa rding cone}.
G\aa rding's beautiful  theory of  hyperbolic  polynomials [G] applies here to give the following.

\Prop{\DD.3.5} {\sl The  G\aa rding cone $\bbm_\PP$ is a convex cone containing the identity
$I$.  It satisfies the property 
$$
\bL_{\PP,k} + \bbm_\PP  \  \ss \  \bL_{\PP,k} \fa k=1,...,m,
$$
that is, $\bbm_\PP$ gives  a monotonicity cone for each of the subequations
$\bL_{\PP,k}$.} In particular, as long as $\bbm_\PP$ contains $\cp$, each branch 
$\L_{Q,k}$ of $Q$ is a subequation.

\medskip
One of the simplest examples comes by taking $\PP(A) = \s_m(A)$, the $m^{\rm th}$ elementary symmetric 
function in the eigenvalues.  Here the G\aa rding cone $\bbm_\PP$ is the set $\{\s_1\geq0,...,\s_m\geq0\}$
(cf. Example \BB.1.2(c)).

In general, for any hyperbolic polynomial $\PP$ as above, one can construct large families of associated
subequations, equipped  with monotonicity cones,  by using  the eigenvalues of $\PP$. For a discussion of this
as well as an elementary introduction to G\aa rding's theory, see [HL$_{7,8}$].


\vskip .3in
\noindent
{\bf \DD.4.  Monotonicity and Duality.} 
The key algebraic fact  that the dual of a translated subequation $F-J$ is just $\wt F -J$
(see \BB.1.5) easily proves the following result, which in turn proves the 
basic algebraic lemmas \DD.1.2 and \DD.2.2.

\Lemma{\DD.4.1} {\sl
Given three subequations $G,M,F \ss J^2(X)$, the fibre-wise sums satisfy:}
$$
G+M\ \ss\ F \qquad\iff\qquad G+\wt F\ \ss\ \wt M.
\eqno{(\DD.4.1)}
$$
\pf 
Note that $J+M\ss F \ \iff\ M\ss -J+F \ \iff\ J+\wt F \ss\wt M$.\qed
 \medskip
 
 Later on, (\DD.4.1) will be implemented  with $G= F^c\ss F$ (cf. (\EE.1.1))
 to obtain weak comparison (see Remark \EE.1.4).


\vfill\eject
\noindent
{\bf \DD.5.  Uniform Ellipticity as Monotonicity.} As noted in Example \DD.4.3 the classical notion
of uniform ellipticity can be reformulated in terms of  monotonicity. We now examine this in greater 
detail.  Suppose that $F$ is a subequation defined on an open set $X\ss \rn$, in the classical way,
  by $F\equiv \{f(x,r,p,A)\geq 0\}$ for a function $f: J^2(X) \to \bbr$ (cf. Appendix A). Then uniform ellipticity 
  (with constants $0<\l < \L$) is the  condition that for  $A,P \in\Symn$ with $P\geq0$,
  $$
  \l \tr (P) \ \leq\  f(x,r,p,A+P) - f(x,r,p, A) \ \leq\ \L\tr( P)
  \eqno{(\DD.5.1)}
  $$
  (and is usually combined with Lipschitz continuity in $p$).  
  This condition can be reformulated in terms of a monotonicity subequation for $F$. To see this 
  it suffices to consider the simplest case $f:\Symn\to\bbr$.  The condition (\DD.5.1)
  is equivalent to requiring that for all $A,B$ (not just $B\geq0$),
  $$
  \cp_{\l,\L}^-(B) \ \leq\ f(A+B) - f(A)\ \leq\   \cp_{\l,\L}^+(B)
 \eqno{(\DD.5.1)'}
  $$
  where $\cp_{\l,\L}^\pm$ are the {\sl Pucci operators} defined by
  $$
  \cp_{\l,\L}^-(B) \ \equiv \l\tr(B^+) + \L\tr(B^-)
  \and
 \cp_{\l,\L}^+ \ \equiv -  \cp_{\l,\L}^-(-B)
  $$
and where $B = B^+ + B^-$ is the decomposition into $B^+\geq0$ and $B^-\leq0$.
It is easy to see that the left hand inequality in (\DD.5.1)$'$ for all $A,B$ is equivalent to the 
right hand inequality for all $A,B$. The desired monotonicity is given by  the {\bf Pucci cone}
$$
\bbp_{\l,\L} \ \equiv \{B\in \Symn :   \cp_{\l,\L}^-(B) \geq0\}.
\eqno{(\DD.5.2)}
  $$
Note that the left hand inequality in (\DD.5.1)$'$ implies the 
monotonicity:
$$
F + \bbp_{\l,\L} \ \ss\ F.
\eqno{(\DD.5.3)}
  $$
The equivalence of $F + \bbp_{\l,\L} \ \ss\ F$ and $\wt F + \bbp_{\l,\L} \ \ss\ \wt F$
corresponds to the equivalence of the right and left hand inequalities in  (\DD.5.1)$'$.

The Pucci cones are convex.  One way to see this is to compute that $\bbp_{\l,\L}$
is the polar of the convex cone on the set  $\{ B\in \Symn : \l I\leq B\leq \L I\}$.

We point out that  Pucci cones provide just one of many choices of a family of 
monotonicity subequations (convex cones) which form a ``fundamental'' neighborhood 
system of $\cp=\{A\geq0\}$, e.g. Example  \DD.3.3  above. All such families give equivalent
notions of uniform ellipticity.

\vfill\eject


\noindent{\headfont \EE.\   Comparison and Strict Approximation.}
\medskip
 Let $F\ss J^2(X)$ be a subequation on a manifold $X$ and for each compact 
 set $K\ss X$ set $F(K) = \USC(K)\cap F(\Int K)$. 
 
 \Def{\EE.1} We say that {\bf comparison holds}
 for $F$ on $X$ if for every compact subset $K$,
the Zero Maximum Principle  
$$
u+v \ \leq\ 0\ \ {\rm on}\ \ \partial K 
\qquad  \Rightarrow\qquad
u+v \ \leq\ 0\ \ {\rm on}\ \  K 
\eqno{(ZMP)}
$$
holds for all
$$
 u\in F(K)\and v\in \ft(K).
 $$

One sees easily that comparison implies {\bf uniqueness for the Dirichlet problem}:
\medskip
\centerline
{
If $u$ and $v$ are $F$-harmonic on $\Int K$ and $u=v$ on 
$\partial K$,
}
\centerline
{ then $u=v$ on $K$}

\vskip .3in
\noindent
{\bf \EE.1.  Weak Comparison.}  A $C^2$ function $u$ on $X$ is said to be 
{\bf strictly} $F$-subharmonic if $J^2_x u\in\Int F_x$ for all $x$.  This notion
has the following  useful extension to functions which are not $C^2$.
For  $c>0$ let $F^c$ be the 
subequation with fibres 
$$
F^c_x \equiv \{J\in F_x : \dist(J, \sim F_x)\geq c\}
\eqno{(\EE.1.1)}
$$
where $\dist$ denotes distance in the fibre $J^2_x(X)$.  This set satisfies conditions (P) and (N).
A function $u\in \USC(X)$ is called {\bf  strictly} $F$-subharmonic  if each $x$ has a neighborhood
$U$ and $c>0$ such that $u$ is $F^c$-subharmonic on $U$.

\Def{\EE.1.1}  We say that {\bf weak comparison holds}
 for $F$ on $X$ if for every compact subset $K$,
$$
u+v \ \leq\ 0\ \ {\rm on}\ \ \partial K 
\qquad  \Rightarrow\qquad
u+v \ \leq\ 0\ \ {\rm on}\ \  K 
$$
holds for all
$$
 u\in F^c(K), \ \ \  v \in \ft(K) \ \ \ {\rm and}\ \ \  c>0.
 $$
We say that {\bf local weak comparison holds}
 for $F$ on $X$ if every point has a neighborhood in which weak comparison holds.
 This weakened form of comparison has several advantages. The first is the following.

 \Theorem{\EE.1.2. (Local implies Global)} {\sl
 If local weak comparison holds  on $X$, then weak comparison holds  on $X$.
 }
 \medskip
 
 A second  important advantage is the following.
 
 \Theorem{\EE.1.3} {\sl
 Suppose $F$ is a subequation on $X$  which is locally jet-equivalent to a constant coefficient
 subequation.  Then weak comparison holds for $F$ on $X$.
 }



\Remark{\EE.1.4} $F^c$ is exactly the subset of $F$ which satisfies the ``weak monotonicity''
$$
F^c +M^c\ \ss\ F \qquad {\rm and\ hence} \qquad F^c+\wt F\ \ss\ {\wt M}^c
$$
where $M^c$ is the universal subequation corresponding to the constant coefficient subequation
$$
\bbm^c\ \equiv\ (-\infty,0]\times \overline{B(0,c)}\times (\cp-c\cdot I).
$$
The  smaller subequation $\bbm_c\ss\bbm^c$
defined by 
$$
\bbm_c \    \equiv \  (-\infty,0]\times \overline{B(0,c)}\times \cp
$$
has dual $\wt{\bbm}_c \supset \wt{\bbm}^c$ which  satisfies the (ZMP). It is the union of three subequations:
$$
\eqalign
{
&\bbr_- \times \rn \times \Symn \qquad{\rm (zeroth\ order)}  \cr
&\bbr \times (\sim{B(0,c)}) \times \Symn \qquad{\rm (dual \ Eikonal)}  \cr
&\bbr \times\rn \times \wt\cp \qquad{\rm (subaffine)},  \cr
}
$$


\vskip .2in
\noindent
{\bf \EE.2.  Strict Approximation.} We say that {\bf strict approximation} holds for $F$ on $X$ if 
for each compact set $K\ss X$, each function $u\in F(K)$ can be uniformly approximated
by functions in $F(K)$ which are strict on $\Int K$.

\Theorem{\EE.2.1} {\sl
If weak comparison and strict approximation hold for $F$ on $X$, then comparison holds 
for $F$ on $X$.
}

\Theorem{\EE.2.2} {\sl
Let $F$ be a subequation on $X$ with a monotonicity cone subequation $M$.
Suppose $X$ carries a $C^2$-function which is strictly $M$-subharmonic.
Then local weak comparison implies global comparison for $F$ on $X$.
}
\medskip
The idea is to approximate $u\in F(K)$ by $u+\e \psi$, $\e>0$,  where $\psi$ is the strictly
$M$-subharmonic function. (The proofs of these theorems can be found in [HL$_6$].)

Thus we see that monotonicity subequations are of central importance in solving the
Dirichlet Problem for nonlinear equations which are degenerate and highly non-convex.

There are times when strict approximation can be achieved by other means. 
One example is given by the Eikonal subequation $|\nabla u|  \leq 1$.
Here the family of functions $u_\e = (1-\e)u$ for $\e>0$ gives strict approximation.


\vskip .3in
\noindent
{\bf \EE.3.  Addition Theorems.} 
In  [HL$_4$]  the following results were proved for pure second-order,
constant coefficient subequations  on an open subset $X\ss\rn$.
We recall that a function $u$ on an open set  in $ \rn$ is {\sl quasi-convex} if
the function $u(x) + c|x|^2$ is convex  for some $c>0$.  Local quasi-convexity is invariant under
coordinate changes and therefore makes sense on manifolds.
\medskip

\noindent

 {\sl
Suppose  $u$ is locally quasi-convex on $X$.  Then}
$$
u \in \bbf(X)\qquad\iff \qquad D_x^2u \in \bbf \ \ {\rm a .e. \ {\sl on} \ \ } X.
$$


\noindent
{\sl
If $\bbf+\bbg   \ss \bbh$,   then for quasi-convex functions $u$ and $v$,}
$$
u \in \bbf(X) \ \ {\sl and}\ \ v\in \bbg(X) \qquad\Rightarrow \qquad u+v\in \bbh(X).
$$

Both of these results hold in much greater generality.

\Theorem{\EE.3.1. (AE Theorem)}  {\sl
Suppose $F$ is a subequation (in the sense of Definition \BB.2.1) on 
a manifold $X$, and suppose  $u$ is locally quasi-convex on $X$.
Then
}
$$
u \in F(X)\qquad\iff \qquad J_x^2u \in F_x   \ \ {\rm a .e. \ on\ \ } X.
$$

\Theorem{\EE.3.2. (Quasi-Convex Addition)}  {\sl
Given three subequations $F$, $G$ and $H$ (as in \EE.3.1) with
$F+G\ss H$, one has that 
$$
u \in F(X) \ \ {\rm and}\ \ v\in G(X) \qquad\Rightarrow \qquad u+v\in H(X).
$$
for locally quasi-convex functions $u$ and $v$.}  

\medskip

Theorem \EE.3.1 follows in an elementary manner from either Jensen's Lemma [J$_1$]
or Slodkowski's Lemma [S$_1$] (in fact, they are equivalent).
Theorem \EE.3.2 is immediate from the first.  These results  will be elaborated
in a forthcoming paper.  

Of course, quasi-convex approximation can be used
in the constant coefficient case to obtain the full Addition Theorem:
$$
u\in \bbf(X) \ \ {\rm and}\ \ v\in \bbg(X) \qquad\Rightarrow \qquad u+v\in \bbh(X).
\eqno{(\EE.3.1)}
$$

\medskip
\noindent
{\bf Application \EE.3.3.  (Comparison via Monotonicity for Constant Coefficient Equations).}
Suppose $\bbf$ satisfies 
$$
\bbf + \bbm \ \ss\ \bbf
\eqno{(\EE.3.2)}
$$
where $\wt \bbm$-subharmonic functions satisfy the Zero Maximum Prinicple.
From  (\EE.3.2) we have  $\bbf+\wt \bbf \ss\wt \bbm$.   Therefore
$$
u \in \bbf(X)\quad{\rm and}\quad v\in\wt\bbf(X)
\qquad\Rightarrow\qquad
u+v\in \wt \bbm(X),
$$
 and so comparison holds for $\bbf$.

Note that  $\bbm$ can be any of the monotonicity cones discussed in Appendix B.
For example, the cone $\bbm= \bbr_-\times \rn\times \cp$ implies comparison 
for all gradient independent subequations.

\vskip.3in



\noindent{\headfont \FF.\  Removable Singularities.}
\medskip
 Monotonicity cones lend themselves nicely to the question of removable singularities 
 for $F$-subharmonic and $F$-harmonic functions.

\vskip .2in
\noindent
{\bf \FF.1.  $M$-Polar Sets.} 
 Suppose $M\ss J^2(X)$ is a convex cone subequation, i.e., one for which the fibres are
 convex cones with vertex at the origin.

 \Def{\FF.1.1}    A closed subset  $E\ss X$  
  is  called $C^\infty$ $M$-{\bf polar} if $E=\{x : \psi(x)=-\infty\}$ for some 
  $M$-subharmonic  function $\psi$ which is smooth on $X-E$.

\medskip\noindent
{\bf Examples.}\smallskip
\noindent (a)  Consider the pure second-order constant coefficient equation $\bbm=\cp$ on $\rn$.
The $\cp$-subharmonic functions are convex (See Proposition \BB.1.7), and so there
do not exist any $C^\infty$ $\cp$-polar sets.

\smallskip
\noindent (b)  Consider the complex analogue  $\cp^\bbc$ on $\bbc^n$.
Then $\cp^\bbc$-subharmonic functions are the standard plurisubharmonic functions
and $\cp^\bbc$-polar sets are standard pluripolar sets. These exist is abundance.
They include, for example,  $\log|f|$ with $f$ holomorphic.

\smallskip
\noindent (c)  For  the quaternionic  analogue  $\cp^\bbh$ on $\bbh^n$
 there is a 2-sphere of complex structures coming from unit 
  imaginary quaternions. A  plurisubharmonic function in any one of these structures
  is $\cp^\bbh$-subharmonic, and so any pluripolar set for that structure is 
$\cp^\bbh$-polar.

\smallskip
\noindent (d) Consider the constant coefficient subequation $\cp(p)$ defined in (\DD.3.1)
and equal to $\bbf(G(p,\rn))$  for integer $p$ (cf. \BB.1.11(c)).
The following result is proved in [HL$_{12}$]  using the theory of  classical Riesz potentials (see  [L] for example).

\Theorem {\FF.1.2}  {\sl Any closed set of locally finite Hausdorff $(p-2)$-measure is $\cp(p)$-polar.}

\vskip .3in
\noindent
{\bf \FF.2.  Removability Results.} The following  removable singularity results on manifolds are proved in
[HL$_{12}$]. Recall  that $M$ is a monotonicity cone for $F$
if and only if it is a monotonicity cone for $\ft$
(see  \DD.2.2).

\Theorem {\FF.2.1} {\sl
Suppose $F$ is a subequation on $X$ with monotonicity cone $M$, and 
 $E\ss X$ is  locally $C^\infty$ $M$-polar with no interior.
Then $E$ is removable for $F$-subharmonic functions which are locally bounded
above across $E$.  More precisely, if  $u\in  F(X-E)$ is locally bounded across $E$, then its
canonical upper semi-continuous extension $U$ to $X$ is $F$-subharmonic on $X$. 
}

\Theorem {\FF.2.2} {\sl
Suppose $F$ is a subequation on $X$ with monotonicity cone $M$, and 
 $E\ss X$ is   locally $C^\infty$ $M$-polar with no interior.
Then   for  $u\in  C(X)$ 
\medskip
\centerline
{
$u$ is $F$-harmonic on $X-E \qquad\Rightarrow\qquad u $ is $F$-harmonic on $X$.
}
}

\medskip

More generally, Theorem \FF.2.1 remains true when $E$ has  interior if the extension $U$ 
is defined to be $\equiv -\infty$ on $\Int E$.

Theorems \FF.2.1 and \FF.2.2 can be applied to the many subequations given in Section \DD.3.
For example, this gives removable singularity results for {\sl all branches} of the homogeneous
complex Monge-Amp\`ere equation  on a complex hermitian manifold.  
 Here $E$ can be any pluripolar
 set (not just a $C^\infty$ pluripolar set). The result also applies to the intrinsic notion of 
 maximal functions on an almost complex manifold (see [HL$_{10}$]).

These general results combined with  Theorem \FF.1.2 above give the following. 
We restrict attention to constant coefficient
pure second-order subequations in $\rn$.
 
 \Cor{\FF.2.3}  {\sl If  $F$ is a subequation for which $\cp(p)$ is a monotonicity cone,
  then any closed set of locally finite Hausdorff $(p-2)$-measure is 
removable for $F$- and $\ft$-subharmonics and $F$-harmonics as in the two theorems above.}

 \medskip
This applies immedately to {\sl all branches} of the equation ${\rm MA}_p$ in Example  \DD.3.2.  It also applies to all  subequations geometrically defined by a subset $\GG$ of the Grassmannian
 $G(p,\rn)$.  (See Example \BB.1.11 and also example (c) following Theorem
\BB.1.12.).  These include 
 the Lagrangian and Special Lagrangian subequations
 in $\bbc^n$, the associative and coassociative subequations in $\bbr^7$,  and the Cayley subequations
 in $\bbr^8$ (where the appropriate value of $p$ is clear in each case).

For the general applicability of this result 
we introduce the following  invariant, which is studied in [HL$_{15}$].

\Def{\FF.2.4}  Suppose $M$ is a convex cone subequation.  The {\bf Riesz characteristic $p_M$
of $M$} is defined to be
$$
p_M\ \equiv\ \sup\{p\in \bbr : I-p P_e \in M \ \forall\ |e|=1\}.
$$

It has the important property that 
$$
\cp(p) \ \ss\ M \qquad\iff\qquad p\ \leq\ p_M.
\eqno{(\FF.2.1)}
$$
and hence:
\medskip
\centerline
 {\sl For any subequation $F$ which is $M$-monotone,  }
\centerline
{ \sl closed sets of locally finite 
Hausdorff $p_M$-measure are $F$-removable as above.}

\medskip
\noindent
{\bf Examples \FF.2.5.}  For  $M = \bbp_{\l,\L}$,  the Pucci cone defined in (\DD.5.2),   the
 Riesz characteristic is
$$
p_M \ =\ {\l\over \L}(n-1) + 1.
$$
 As a consequence one retrieves the removable singularity results in [AGV]. In fact Corollary
 \FF.2.3 is stronger since it applies to interesting equations which are not uniformly elliptic.
 
 For $M=\cp(\d)$, another choice for defining uniform ellipticity, the Riesz characteristic is
 $$
 p_M\ =\ {\d n +1\over \d+1}
 $$

\medskip
\noindent
{\bf Final Remark.}  In the special case of convex subequations (in the general setting of manifolds) 
there are many 
interesting  removability results [HL$_{12}$].  They come from combining the   Strong Bellman Principle
(see \S 10) and known results ([Le], [HP$_{1,2}$], [H], [S]) for linear elliptic equations. 
See   [HL$_{13}$] for details.

\vfill\eject


\noindent{\headfont \GGG.\  Boundary Convexity.}
\medskip

 Fix a subequation $F$ on a manifold $X$ and   a domain $\O\ss\ss X$ 
 with smooth boundary. We shall be interested in the Dirichlet problem
 for $F$-harmonic functions on $\O$.
In this chapter we present geometric conditions on $\bo$ 
which guarantee  the existence of solutions for all continuous
boundary functions.
These conditions are based on the following concept.

\vskip .3in
\noindent
{\bf \GGG.1.  The Asymptotic Interior of a Reduced Subequation.} 
Throughout this section we assume that $F$ is a  subequation which is ``independent
of the $r$-variable''  or ``reduced''.  This means that  with respect to 
the splitting 

\smallskip
\centerline
{$J^2(X) = \bbr\oplus J^2_{\rm red}(X)$
}
\smallskip
\noindent
 in (\CC.1.1), $F$ is of the form
$F=\bbr\times F_0$.  For simplicity we just take $F\ss J^2_{\rm red}(X)$.

 \Def{\GGG.1.1}   
 The {\bf asymptotic interior} $\Fa$ of $F$ is the set of all $J\in J^2_{\rm  red}(X)$ 
 for which there exists a neighborhood $\cn(J)$ in the total space of $J^2_{\rm  red}(X)$
 and a number $t_0> 0$ such that 
 $$
 t\cdot \cn(J)\ \ss\ F\fa t\geq t_0
 $$
\smallskip

The set $\Fa$ is an open cone in $J^2_{\rm red}(X)$ 
which satisfies Condition (P).
If $F$ is itself a cone, then $ \Fa\ =\ \Int F$.  Otherwise,  $\Fa$ is smaller than $\Int F$  
and may be empty.

 \Def{\GGG.1.2}   
A function $u\in C^2(X)$  is called {\bf strictly $\Fa$-subharmonic} if $J^2_{{\rm red}, x} u \in \Fa$ for all $x$. 

\medskip

Let $\O \ss X$ be a domain with smooth boundary $\bo$.  By a {\sl defining function} for 
$\bo$ we mean a smooth function $\rho$ defined on a neighborhood of $\bo$ such that
$\bo = \{x:\rho(x)=0\}$, $d\rho\neq 0$ on $\bo$, and $\rho<0$ on $\O$.

 \Def{\GGG.1.3}   Suppose  $F$ is a reduced subequation.  The boundary $\bo$ is said to be 
 {\bf strictly $F$-convex at}  $x\in\bo$ if  there exists a  strictly $\Fa$-subharmonic
 defining   function for $\bo$ on some neighborhood of $x$

\medskip
This is equivalent to either of the following two conditions.
\smallskip\noindent (i)
For   some local defining function $\rho$,\ \ 
 $J^2_{{\rm red}, x} \rho  \in \Fa$. 
 
 \smallskip\noindent (ii)
For   any  local defining function $\rho$, \ \ 
 $J^2_{{\rm red}, x} \rho + t(d\rho)_x\circ (d\rho)_x  \in \Fa$ for all $t\geq $ some $t_0$.

\vskip .3in
\noindent
{\bf \GGG.2.  General $F$-Convexity.} 
Suppose now that $F\ss J^2(X)$ is a general subequation on $X$.  For each $\l\in \bbr$  there is a reduced
subequation $F_\l \ss J^2_{\rm red}(X)$ obtained by fixing the $r$-variable to be $\l$,
that is
$$
F_\l \ \equiv \ F\cap \left ( \{\l\}\times J^2_{\rm red}(X)\right).
$$

As above we fix a domain  $\O\ss X$  with smooth boundary $\bo$.

\Def{\GGG.2.1} Suppose  $F$ is a general  subequation.  The boundary $\bo$ is said to be 
 {\bf strictly $F$-convex at}  $x\in\bo$ if it is strictly $\oa{{F_\l}}$-convex at $x$ for all $\l\in \bbr$.

\medskip
For example, consider the universal riemannian subequation $F$
given by $\Hess \,  u \geq0$ and $\det\{\Hess\, u\} \geq e^u$.  
Then $F_\l$ is given by the condition that $\Hess \,  u \geq0$ and
$\det\{\Hess\, u\} \geq e^\l$.
 One easily checks that for every $\l$, $\oa{{F_\l}}$ is the open cone  
 $
\{ \Hess\ u>0\},
 $
and so in this case the strictly $F$-convex boundaries are just the classical strictly convex boundaries.

Strict $F$- and $\ft$-convexity of $\bo$  at each point are sufficient  for the construction
of barriers used in the proof of the existence of solutions to the Dirichlet problem.

\vskip .3in
\noindent
{\bf \GGG.3.  $F$-Convexity in Terms of the Second Fundamental Form.} 
For a reduced subequation $F$ on a riemannian manifold $X$, 
 the $F$-convexity of a boundary $\bo$  can be characterized in terms of 
 its  second fundamental form $II_{\bo}$ with respect to  the outward-pointing unit  normal $\nu$. 
We use the decomposition given by (\BB.2.4): 
  $$
  J^2_{\rm red}(X) = T^*X\oplus
 \Sym(T^*X).
 $$

 \Prop {\GGG.3.1}  {\sl  The boundary $\bo$ is strictly $F$-convex at $x\in \bo$
 if and only if   
 $$  
 \left(  \nu, \,  tP_\nu \oplus II_{\bo}\right) \ \in\ \oa{{F_x}} \ \ \ 
 \fa t  \geq {\rm some\ } t_0.
 \eqno{(\GGG.3.1)}
 $$
 where $P_\nu$ denotes orthogonal projection onto the normal line $\bbr \nu$ at $x$.}
\medskip\noindent
{\bf Note.} 
Blocking   with respect to the decomposition
 $T_xX = \bbr \nu \oplus T_x(\bo)$, (\GGG.3.1) can be rewritten
  $$  
 \left(  (1,0), \left(\matrix {t & 0 \cr 0 & II_{\bo}\cr}\right)\right) \ \in\      \oa{{F_x}} \ \ \ \ 
 \fa t \geq {\rm some\ } t_0.
 \eqno{(\GGG.3.2)}
 $$

\vskip .3in
\noindent
{\bf \GGG.4.  Examples.} 
\medskip

(a)\ \ {\bf $k$-Laplacians.}
There are many  examples where every boundary is strictly $F$-convex.
The simplest one is the subequation $\D u\geq 0$ or more generally
$\D u \geq f(x,u)$ where $f$ is non-decreasing in $u$.  

Other examples come from the constant coefficient 
 $k$-Laplace subequation,  defined by 
$$
\bbf_k^{\rm Lap} \ \equiv \ {\rm Closure} \left \{(p,A) : |p|^2 \tr\, A + (k-2)\,p^t A p> 0\right \}
 \eqno{(\GGG.3.3)}
 $$  
 where $k\geq 1$. These equations are self-dual.
 Since $\bbf_k^{{\rm Lap}}$ is a cone, $\oa{\bbf_k}^{{\rm Lap}} = \Int \, \bbf_k^{{\rm Lap}} $.
 One can check directly from (\GGG.3.2) that for $k>1$ every boundary is $\oa{\bbf_k}^{{\rm Lap}}$-convex.

When $k=1$ this equation is the implicit minimal surface equation  studied by De Giorgi and his school
[Giu]. Here one sees that a boundary $\bo$ is strictly $\bbf_1^{\rm Lap}$-convex if and only if it is
strictly mean convex, i.e., $\tr (II_{\bo} ) >0$ at all points.

At the other extreme is the infinity Laplacian (cf.  [CIL], [J$_2$], [ESm])
$$
\bbf_{\infty}^{\rm Lap} \ \equiv \ {\rm Closure} \left \{(p,A) : p^t A p> 0\right \}
 \eqno{(\GGG.3.4)}
 $$  
where again all boundaries are strictly $\bbf_{\infty}^{\rm Lap} $-convex.

\medskip
(b)\ \ {\bf Elementary Symmetric Functions of Hess$(u)$}.
Consider  Example \BB.1.2(c)
$$
\bbf_{\s_k} \ \equiv\ \{\s_k(A)\geq0, \ \s_{k-1}(A)\geq0,\  ... \ , \s_1(A)\geq0\}
 \eqno{(\GGG.3.5)}
 $$ 
 which can be extended to the complex and
quaternionic cases, and carried over to riemannian manifolds.
One finds that $\bo$ is strictly $\bbf_{\s_k}$-convex if and only if 
$$
\s_{k-1}\left( II_{\bo} \right) > 0,\ \ \s_{k-2}\left( II_{\bo} \right) > 0, \ ... \ , \ \s_{1}\left( II_{\bo} \right) > 0.
$$
Moreover, if $\bo$ is strictly $\bbf_{\s_k}$-convex, then it is $\bbf_{\s_k,i}$-convex
for every branch $\bbf_{\s_k,i}$ of the equation $\s_k(\Hess\,u)=0$ (see Section \DD.3).
This includes the dual subequation $\wt\bbf_{\s_k}$, which is the bottom branch.

\medskip
(c)\ \ {\bf Geometrically Defined Subequations}.
Consider now the subequations discussed in Example \BB.1.11.
Here the boundary convexity is particularly nice.  Fix a compact subset  $\GG\ss G(p,\rn)$ 
 and define $\bbf(\GG)$ as in (\BB.1.6).  Then a boundary $\bo$ is strictly
$\bbf(\GG)$-convex if and only if 
$$
\tr_W\left\{  II_{\bo}\right\} \ > \ 0
\fa \GG\, {\rm planes} \ W \ {\rm which\ are\  tangent \ to\ }\bo.
 \eqno{(\GGG.3.6)}
 $$ 
This condition holds automatically  at  $x\in\bo$ if there 
are no $\GG$-planes tangent to $\bo$ at $x$. 

On the other hand, if $\GG=G(p,\rn)$, then $\bo$ is strictly $F(\GG)$-convex if and only if 
$II_{\bo}$ has positive trace on all tangent $p$-planes,  i.e., $\bo$ is {\sl $p$-convex}  as in  [Wu], [Sha$_{1,2}$].

For example, suppose $\GG\ss G(1,\bbr^2)$ is the single point $\GG =\{x$-axis$\}$.
Then a domain $\O\ss\ss\bbr^2$  with smooth boundary is strictly $\GG$-convex iff
the curvature vector of $\bo$ points strictly inward at every horizontal tangent.
This implies that all horizontal slices of $\O$ are connected.  Thus, one can see directly that  the 
Dirichlet problem for $\GG$-harmonic functions ($u_{xx}=0$) is uniquely solvable 
for all continuous boundary data. 

A classical example comes from the set $\GG=G_\bbc(1,\bbc^n)\ss G(2,\bbr^{2n})$ of complex lines
in $\bbc^n$. A domain $\O\ss \bbc^n$ is strictly $\GG$-convex iff it is strictly pseudo-convex
in the usual sense in complex analysis (cf. [Ho$_1$]). This is the boundary convexity
required to solve the Dirichlet problem for $\cp^\bbc = \bbf(\GG)$-harmonic functions,
i.e., for solutions to  the homogeneous complex Monge-Amp\`ere equation.

We note that in all cases $F(\GG)\ss\wt{F(\GG)}$, so that a strictly $F(\GG)$-convex boundary
is automatically  strictly $\wt{F(\GG)}$-convex.

\medskip
(d)\ \ {\bf p-Plurisubharmonic Functions}.
Consider now the $p^{\rm th}$ branch of the homogeneous complex Monge-Amp\`ere
equation.  This is the pure second-order subequation given by
$\L^\bbc_{p} \equiv \{A:  \l^\bbc_{p}(A)\geq 0\}$
where $\l^\bbc_{1}(A)\leq \cdots\leq \l^\bbc_{n}(A)$ are the ordered eigenvalues of the 
hermitian symmetric part of $A$ (see \BB.1.3 and \BB.1.10). 
The $\L^\bbc_{p}$-subharmonic functions are the classical $(p-1)$-plurisubharmonic functions
in complex analysis --  those for which the complex hessian has at least $n-p+1$ non-negative eigenvalues.
The Dirichlet problem for $\L^\bbc_{p}$-harmonic functions was studied by Hunt and Murray [HM] 
and then solved by Slodkowski [S$_1$]. A smooth boundary $\bo\ss\bbc^n$ is strictly $\L^\bbc_{p}$-convex iff 
$$
\l_p^\bbc\left( II_{\bo}\right) \ \geq \ 0,\qquad{\rm or\ equivalently}
 \eqno{(\GGG.3.7)}
 $$ 
\centerline
{the Levi form of $\bo$ has $n-p-1$ eigenvalues $\geq0$ 
at each point.}

\medskip
(e)\ \ {\bf Calabi-Yau-Type Equations.}
Let $X$ be a complex hermitian manifold.  
Consider the  subequation $F$ on $X$ corresponding to  $\det_\bbc(I+\Hess_\bbc u) \geq f(x,u)$ 
for a continuous  $f > 0$ which is non-decreasing in $u$ and $I+\Hess \, u \geq0$. 
For $\l\in\bbr$ the subequation $F_\l$ given in Section \GGG.2 corresponds to 
$\det_\bbc(I+\Hess_\bbc u) \geq f(x,\l)$ at each point. One checks that $F_\l$-convexity
of  a boundary  $\bo$ amounts to the statement that $(II_{\bo})_\bbc > -I$ at each point
(a condition independent of $\l$). Levi convexity of the boundary ($(II_{\bo})_\bbc > 0$)
will certainly suffice.

\medskip
(f)\ \ {\bf Principal curvatures of the graph.} Other equations of interest are those which impose
conditions on the principal curvatures of the graph of the function $u$ in $X\times \bbr$.
See [HL$_6$, \S 11.5] for a complete discussion of this case.

\vfill\eject


\noindent{\headfont \HH.\  The Dirichlet Problem.}
\medskip
 Throughout this chapter $F\ss J^2(X)$ will be a subequation on a manifold 
 $X$ and $\O\ss\ss X$ will be a domain with smooth boundary $\bo$.
 We shall say that {\sl existence holds for the Dirichlet Problem}
 for $F$-harmonic functions on $\O$ if for each continuous function $\vf\in C(\bo)$
 there exists a function $u\in C(\ob)$ such that
 \medskip
\qquad\qquad\qquad\qquad (i)\ \ $u$ is $F$-harmonic on $\O$, and
 \medskip
\qquad\qquad \qquad\qquad (ii)\ $u\bigr|_{\bo} = \vf$.\bigskip 
\noindent
We say that {\sl uniqueness holds for this problem} if for each $\vf\in C(\bo)$, 
there exists at most one such function $u$.

\vskip .3in
\noindent
{\bf \HH.1.  General Theorems.} 
It is an elementary fact that if comparison holds for $F$ on $X$ (see Definition \EE.1),
then uniqueness holds for the Dirichlet problem. Under appropriate boundary convexity
comparison also implies existence.

\Theorem{\HH.1.1} {\sl   Suppose comparison holds for $F$ on $X$.  Then existence and uniqueness
hold for the Dirichlet problem for $F$-harmonic functions on any domain $\O\ss\ss X$ whose boundary
is both strictly $F$-convex and strictly $\ft$-convex.
}

\medskip

Note that $u$ is $F$-harmonic  if and only if $- u$ is $\ft$-harmonic. Thus, it is expected that 
both conditions,  strict $F$ and $\ft$ convexity, should be required, if one  of them is.
Often  one of these convexity conditions implies the other. This is clearly the case
for $F=\cp$ in $\rn$ where strict $\cp$-convexity is the usual strict convexity and $\wt\cp$-convexity
is much weaker.  It also holds  in the case of $q$-\psh functions  (Example \GGG.4(d))
where by (\GGG.3.7) $\cp^\bbc_q$-convexity implies $\cp^\bbc_{q'}$-convexity if $q < q'$.
 This is reflected in the work of Hunt and Murray [HM] who noted the failure of the statement
 when only one convexity condition is required.

Theorems \EE.1.2 and \EE.2.1 imply that
\medskip
\centerline
{
\sl If local weak comparison and strict approximation hold for 
$F$ on $X$, }
\centerline{\sl
then comparison holds for $F$ on $X$.}

\Theorem{\HH.1.2} {\sl
Let $F$ be a subequation with monotonicity cone $M$.  Suppose that:
\medskip

(i)\ \   $F$ is locally affinely  jet-equivalent to a constant coefficient subequation, and

\medskip

(ii)\ \ $X$ carries a strictly $M$-subharmonic function.
\medskip\noindent
Then existence and uniqueness
hold for the Dirichlet problem for $F$-harmonic functions on any domain
$\O\ss\ss X$  whose boundary
is both strictly $F$-  and   $\ft$-convex.
}

\medskip
Comparison and therefore uniqueness follow from Theorems \EE.1.3 and \EE.2.2.
It is then proved,  using comparison and barriers constructed from boundary convexity, that 
existence also holds.  Further details are given in \S 8.
\medskip

Assumption (ii) is always true  for pure second-order equations in $\rn$ 
(and in any complete simply-connected  manifold of non-positive sectional
curvature) since the subequation $\cp$ is always a monotonicity cone by the positivity condition (P)
and $|x|^2$ is strictly $\cp$-convex.

On the other hand something like assumption (ii) must be required in the general case.
For example, suppose $F$ is a universal riemannian equation as in \BB.2.3.
One could completely change the geometry (and topology) of the interior of a domain
$\O\ss X$ without changing the $F$-convexity of the boundary.  Take the subequation
$\cp$ on the euclidean ball, and change the interior so that it is not contractible.  Then 
there are no $\cp$-subharmonic (riemannian convex) functions on the resulting space,
and certainly no $\cp$-harmonic ones.

In homogeneous spaces one can apply a trick of Walsh [W] to establish 
existence without uniqueness.

\Theorem{\HH.1.3} {\sl
Let $X=G/H$ be a riemannian homogeneous space and suppose that $F\ss J^2(X)$ is a 
subequation which is invariant under the natural action of $G$ on $J^2(X)$.
Let $\O\ss\ss X$ be a connected domain whose boundary is both $F$ and $\ft$
strictly convex. Then existence holds for the Dirichlet problem for $F$-harmonic functions on $\O$.
}

\medskip
This theorem applies to give (the known) existence for the $k$-Laplacian, $1<k\leq \infty$ on arbitrary domains,
and for the 1-Laplacian on mean convex domains in  $G/H$.
The literature on these equations in $\rn$ is vast.
See [JLM], [CIL], [J$_2$], [ESm] and references therein, for example.  
We note that even in $\rn$, uniqueness for the 1-Laplacian  fails catastrophically.  For a generic smooth
function on the boundary of the unit disk in $\bbr^2$ there are families of solutions
to the Dirichlet problem parameterized by $\bbr$ (and often $\bbr^m$ for large $m$)!

The proof of existence in the theorems above uses the standard Perron method
based on the properties in Theorem \BB.3.1.
 Given $\vf\in C(\bo)$, consider the family
 $$
 \cf(\vf) \ \equiv \ \{u\in\USC(\ob)\cap F(\O) : u\leq \vf \ \ {\rm on}\ \bo\},
 $$
and define the {\sl Perron function} to be the upper envelope of this family:
$$
U(x)\ \equiv\ \sup_{u\in \cf(\vf)} u(x).
\eqno{(\HH.1.1)}
$$

\Prop{\HH.1.4} {\sl Suppose  that $F$ satisfies weak comparison and that $\bo$ is 
both $F$ and $\ft$ strictly convex.  Then the upper and lower semi-continuous
regularizations $U^*$ and $U_*$ of $U$ on $\ob$ satisfy:
\medskip
\hskip.4in (i) \ \ $U^*=U_* = U = \vf$  on $\bo$,

\medskip
\hskip.4in (ii) \ \ $U=U^*$ on $\ob$

\medskip
\hskip.4in (iii) \ \ $U$ is $F$-subharmonic and $-U_*$ is $\ft$-subharmonic on $\O$.
}
\medskip

The classical barrier argument, used by  Bremermann [B] for the case $F=\cp^\bbc$,
establishes (i), while weak comparison is used in (ii). 
Part (iii) relies on a ``bump argument'' found in Bedford and Taylor [BT$_1$] and also in [I].

When one can ultimately establish comparison, as in Theorem \HH.1.2,
the Perron function is the unique solution.  When this is not necessarily
possible, as in Theorem \HH.1.3, arguments of Walsh [W] can be applied
to show that the Perron function is a solution.

In this latter case one can say more.  Fix $F$  and $\O$ as in Theorem \HH.1.3.
Suppose
\smallskip
\centerline
{
$U$ is the Perron function for $F$ on $\O$ with boundary values $\vf$, and
}

\smallskip
\centerline
{
$-\wt U$ is the Perron function for $\wt F$ on $\O$ with boundary values $-\vf$. \ \ \ \ \ \ 
}
\smallskip
\noindent
Both $U$ and $\wt U$ solve the Dirichlet problem for $F$-harmonic functions on $\O$
with boundary values $\vf$, and if $u$ is any other such solution,
$$
\wt U\ \leq \ u\ \leq \ U.
\eqno{(\HH.1.2)}
$$

Theorems \HH.1.2 and \HH.1.3 have wide applications.
In the following sections we will examine some specific examples.


\vskip .3in
\noindent
{\bf \HH.2.  Manifolds with Reduced Structure Group.}  Fix a constant coefficient subequation
$\bbf\ss\bbj^2$, and let
$$
G\ \equiv\ G_\bbf\ \equiv\ \{g\in {\rm O}(n) : g(\bbf)=\bbf\}
\eqno{(\HH.2.1)}
$$
where $g$ acts naturally on $\bbj^2$ by $g(r,p,A) = (r,gp, g^t Ag)$.

\Def{\HH.2.1} Let $X$ be a riemannian $n$-manifold and $G\ss {\rm O}(n)$ a subgroup.
A {\bf topological $G$-structure} on $X$ is a family $\{(U_\a, e_\a)\}_\a$ where 
$\{U_\a\}_\a$ is an open covering of $X$ and  each  $e_\a = (e_\a^1,...,e_\a^n)$ is a
continuous   
tangent frame field on $U_\a$, such that for all $\a,\b$ the change of framing
$g:U_\a\cap U_\b \to {\rm O}(n)$ takes values in $G$.
\medskip

Each  constant coefficient subequation  $\bbf$ canonically determines
a subequation $F$ on any riemannian manifold $X$ equipped with  a topological
$G_\bbf$-structure. (Use the splitting (\BB.2.4) and then the  
trivializations induced  by the local tangent frames. The subequation 
determined by $\bbf$ in these trivializations is preserved under the change of framings.) 
By Proposition \CC.2.5, $F$ is locally jet-equivalent to $\bbf$.

If $\bbm$ is a $G_\bbf$-invariant monotonicity cone for $\bbf$, then
the  corresponding subequation $M$ on $X$ is a monotonicity cone for $F$.
Note that the maximal monotonicity cone for $\bbf$ is always $G_\bbf$-invariant.

\Theorem{\HH.2.2}  {\sl
Let $F$ be a subequation with monotonicity cone $M$ canonically determined
by 
$\bbf$ and $\bbm$ on a riemannian manifold $X$ with a topological
$G_\bbf$-structure. Let $\O\ss\ss X$ be a domain with smooth boundary which is both
$F$ and $\ft$ srictly convex. Assume there exists a strictly $M$-subharmonic function
on $\ob$.  Then the Dirichlet Problem for $F$-harmonic functions is uniquely solvable 
for all  $\vf\in C(\bo)$.
}
 
\Ex{\HH.2.3}  

\medskip
\noindent
{\bf (a)  Universal Riemannian Subequations:}
As noted in Remark \BB.2.3, 
if    $G_\bbf ={\rm O}(n)$, then $\bbf$ universally determines a subequation on
every riemannian manifold by choosing the framings  $e_\a$ to be orthonormal.
In particular this covers all branches of the homogeneous Monge-Amp\`ere equation.
In fact, it covers all pure second-order subequations which depend only on the 
ordered eigenvalues of the Hessian. The subequation $\cp= \{\Hess\,u \geq0\}$
is a monotonicity cone for all such equations. Thus Theorem \HH.2.2 applies  to all such $F$'s
in any region of $X$ where there exists a smooth strictly convex function.

Other interesting examples are given by the branches of the $p$-convex Monge-Amp\`ere 
equation ${\rm MA}_{p}$
given in example \DD.3.2. Here  the monotonicity cone is $\cp(p)$, and the appropriate 
boundary convexity is the $p$-convexity discussed in \GGG.4 (c).

Further examples come from elementary symmetric functions of $\Hess \,u$ (see 
\GGG.4 (b) and the discussion after \DD.3.5.), and functions of eigenvalues of the
graph (\GGG.4 (f)).

\medskip
\noindent
{\bf (b)  Universal Hermitian Subequations:} If $G_\bbf  =  {\rm U}(n)$, then $\bbf$
universally determines a subequation on  every almost complex
hermitian manifold.   For example, this covers all pure second-order subequations 
which depend only on the  ordered eigenvalues of the hermitian symmetric part 
$\Hess_\bbc u$ of  $\Hess\,u$. For such equations, $\cp^\bbc = \{\Hess_\bbc u\geq0\}$ is 
a monotonicity cone.  Thus, for example, one has the following consequence of  Theorem
\HH.2.2.  {\sl Let $X$ be an almost complex hermitian manifold, and $\O\ss\ss X$ a 
smoothly bounded domain with a strictly plurisubharmonic ($\cp^\bbc$-subharmonic)
defining function.  Then the Dirichlet problem for every branch of the homogeneous 
complex Monge-Amp\`ere equation is uniquely solvable on $\O$.}

A similar result holds for branches of the equation  ${\rm MA}_p^\bbc$ where $p$-convexity 
of the Levi form on the boundary plays a role (see \GGG.4 (d)).

The discussion of elementary symmetric functions also carries over to this case,

Theorem \HH.2.2 can similarly be applied to Calabi-Yau type equations (\GGG.4 (e)).

All of this discussion can be replicated for almost quaternionic hermitian manifolds.

\medskip
\noindent
{\bf (c)  Geometrically Defined Subequations:} Theorem \HH.2.2 applies directly to all
subequations geometrically defined by a compact subset $\GG\ss G(p,\rn)$
(see  \BB.1.11, \BB.1.12 and  \GGG.4 (b)). {\sl Suppose $X$ has a topological $G$-structure
where $G=\{g\in {\rm O}(n) : g(\GG)=\GG\}$ and let $F(\GG)$ be the corresponding
subequation on $X$.  Suppose $\O\ss X$ is a domain with a 
{\sl global}  defining function which is strictly $\GG$-\psh.
Then the Dirichlet problem for $\GG$-harmonic functions is uniquely solvable on $\O$.}

Thus, one can solve the Dirichlet problem for (in fact, all branches of) the Lagrangian harmonic equation
(see \BB.1.11 (d)) on domains with a strictly Lagrangian-\psh defining function.

One can also solve for $G(\vf)$-harmonic functions  on strictly $G(\vf)$-convex domains
in a manifold with a topological calibration $\vf$.
A typical example is the following. Let $X$  be a riemannian 7-manifold with a topological
G$_2$-structure determined by a global associative 3-form $\vf$ of constant comass 1
(Such structures exist on $X$ if and only if $X$ is a spin manifold.)
Then the Dirichlet problem for $G(\vf)$-harmonic functions is uniquely solvable on any
domain with a strictly $G(\vf)$-\psh defining function.


\vskip .3in
\noindent
{\bf \HH.3.  Inhomogeneous Equations.}     Since  Theorem \HH.1.2 assumes {\sl affine}
jet-equivalence, it applies to inhomogeneous equations as in Examples \CC.2.7-8.
In these cases boundary convexity and monotonicity cones are the same as in the 
homogeneous case.


\vskip .3in
\noindent
{\bf \HH.4.  Existence Without Uniqueness.}  Theorem \HH.1.3 applies in   cases where
monotonicity cones do  not exist, such as the $1$-laplacians in \GGG.4 (a). 
As previously noted, solutions of the Dirichlet problem
 for the 1-laplacian are highly non-unique. However, they are all caught between the
 Perron functions $U$ and $\wt U$ (see (\HH.1.2) above).


\vskip .3in
\noindent
{\bf \HH.5.  Parabolic Equations.}  
The methods and results above carry over effectively to parabolic equations.
Let  $X$ be a riemannian $n$-manifold  with a topological 
$G$-structure  for $G\ss {\rm O}(n)$, and consider a constant coefficient subequation of
the form
$$
\bbf \ =\ \{J\in \bbj^2 : f(J)\geq0\}
$$
where $f: J^2(X) \to \bbr$  is $G$-invariant,  $\cp$- and $\cn$-monotone,
and Lipschitz in the reduced variables $(p,A)$. 
This induces a subequation $F$ on $X$.
The associated constant coefficient parabolic subequation ${\bf H}_{\bf F}$ on $\bbr\times \rn$
is defined by
$$
f(J)-p_0\ \geq\ 0
$$
(where $p_0$ denotes the $u_t$ component of the 2-jet of $u$),
and it induces the {\sl associated parabolic subequation} $H_F$
 on the riemannian product $\bbr \times X$.
The $H_F$-harmonic functions
are solutions of the equation 
$$
u_t = f(u, Du, D^2 u).
$$ 
Examples which can be
treated include:

\smallskip
 (i) \  $f = \tr A$, the standard heat equation $u_t = \D u$ for  the Laplace-Beltrami operator.

\smallskip
 (ii)     $f = \l_q(A)$, the $q$th ordered eigenvalue of $A$.  This is 
 the natural parabolic equation associated to the $q$th branch of the Monge-Amp\`ere equation.

\smallskip
 (iii)  \  $f = \tr A +  { k\over  |p|^2+\e^2}  p^t A p$ for $k\geq - 1$ and $\e > 0$.
When $X=\rn$ and $k=-1$, 
the solutions $u(x,t)$ of the associated parabolic equation, in the limit as
$\e\to0$, 
have the property that the associated level sets
$
\Sigma_t\ \equiv\ \{x\in\rn : u(x,t)=0\}
$
are evolving by mean curvature flow (cf.  [ES$_*$],  [CGG$_*$],  [E] and [Gi].)

\smallskip
 (iv)  $f = \tr\{\arctan A\}$.  
 When $X=\rn$,  solutions $u(x,t)$ 
 have the property that the graphs of the gradients:
 $
 \G_t\ \equiv\ \{(x,y)\in \rn\times \rn =\bbc^n : y = D_xu(x,t)\}
 $
are Lagrangian submanifolds which evolve the initial data by 
mean curvature flow.  (See [CCH].)
\smallskip

Techniques discussed above  show  that:\medskip
\centerline
{\sl Comparison holds for the subequation $H_F$ on $X\times \bbr$.
}
\medskip

Applying standard viscosity techniques for parabolic equations, one can prove more. 
Consider a compact subset  $K \ss \{t\leq T\}\ss X\times \bbr$ and 
let $K_T \equiv K\cap \{ t=T\}$ denote the terminal time slice of $K$.  Let
$
\partial_0 K \  \equiv \ \partial K  - \Int K_T
$
denote the {\sl parabolic boundary of $K$}.  Here $\Int K$ denotes the relative interior in 
$\{t=T\}\ss X\times\bbr$.
 We say that {\sl parabolic comparison holds for $H_F$} if for all such $K$ (and $T$)
$$
u+v\ \leq\ c \quad{\rm on}\ \ \partial_0 K\qquad\Rightarrow\qquad u+v\ \leq\ c \quad{\rm on}\ \ \Int K
$$
for all $u\in H_F(K)$ and $v\in \wt H_F (K)$.  Then one has that:
\medskip
\centerline{\sl Parabolic comparison holds for  the subequation $H_F$ on $X\times \bbr$.
}
\medskip

Under further mild assumptions on $f$ which are satisfied in the examples above, 
one also has existence results.   Consider a domain $\O \ss X$ whose boundary is
strictly $F$- and ${\wt F}$-convex.  Set $K= \overline\O \times [0,T]$.
Then 
\medskip
\centerline{\sl For each $\vf\in C(\partial_0 K)$ there exists a unique function $u\in C(K)$ such
that
}
\centerline{\sl $u\bigr|_{\Int K}$ is $H_F$-harmonic\ \  and\ \  $u\bigr|_{\partial_0 K} = \vf$.
}
\medskip
\noindent
One then obtains corresponding long-time existence results.


\vskip .3in
\noindent
{\bf \HH.6.  Obstacle Problems.}  
The  methods discussed here lend themselves easily to solving boundary
value problems with obstacles.  Suppose that $F=\bbr\times F_0$ is a reduced subequation,
i.e., independent of the $r$-variable.  Given $g\in C(X)$, the {\sl  associated obstacle subequation}
is defined to be 

\smallskip

\centerline
{$H\equiv (\bbr_-+g)\times F_0$\ \ \  where $\bbr_- \equiv \{r\leq0\} \ss\bbr$.
}
\smallskip
The following facts are easy to prove.
\smallskip

\noindent
$\bullet$ \ \ The $H$-subharmonic functions  are  the $F$-subharmonic functions
$u$ which satisfy $u\leq g$.

\smallskip

\noindent
$\bullet$ \ \ If $F$ has a monotonicity cone $M =  \bbr\times  M_0$, then 
$M_-\equiv \bbr_-\times M_0$ is a monotonicity cone for $H$.

\smallskip

\noindent
$\bullet$ \ \ If  $X$ carries a strictly $M$-subharmonic function $\psi$, then on any given compact
set, the function  $\psi-c$ is strictly $(M_-)$-subharmonic for $c>0$ sufficiently large.

\smallskip

\noindent
$\bullet$ \ \ If  $F$ is locally affinely jet-equivalent  to a constant coefficient reduced subequation
$\bbr\times \bbf_0$, then $H$ is locally affinely jet-equivalent to the subequation $\bbr_-\times \bbf_0$.

\smallskip

Consequently, under the assumptions in Theorem  \HH.1.2 on a
 reduced subequation $F=\bbr\times F_0$ with monotonicity  cone $M=\bbr\times M_0$, 
{\bf  comparison holds for each associated obstacle subequation $H\equiv (\bbr_-+g)\times F_0$.}

However, existence fails for a boundary function $\vf\in C(\bo)$ unless 
 $\vf \leq g\bigr|_{\bo}$.  Nevertheless, {\bf   if $\bo$ is both $F$ and $\wt F$ strictly convex as in
 Theorem \HH.1.2, then existence holds for each boundary function $\vf \leq g\bigr|_{\bo}$. }

 To see  that this is true, note the following. The Perron family for a boundary function
 $\vf\in C(\bo)$ consists of those $F$-subharmonic functions $u$ on $\O$ with 
$u\bigr|_{\bo} \leq \vf$ (the usual family for $F$) subject to the additional constraint
$u\leq g$ on $\O$. The dual subequation to $H$ is 
$\wt H = [(\bbr_--g)\times J^2_{\rm red}(X) ] \cup \ft$ so that the boundary $\bo$ is strictly
$\wt H$-convex if  it is strictly $\ft$-convex.  Although $\bo$ can never be strictly
$H$-convex (since $(\oa{{F_\l}})_x = \emptyset$ for $\l> g(x)$), the only place that this hypothesis is used 
in proving  Theorem \HH.1.2 for $H$ is in the barrier construction which appears in the proof of 
 Proposition $F$ in [HL$_6$].
However, if $\vf(x_0) \leq g(x_0)$, then the barrier $\b(x)$ as defined in (12.1) in [HL$_6$]
 is not only $F$-strict near $x_0$ but also automatically $H$-strict since $\b < g$.

The obstacle problem for the basic subequation $\cp$ is related to convex envelopes.
This was discovered by Oberman [O] and developed by Oberman-Silvestre  [OS].


  \vskip .3in


\noindent{\headfont \II.\  Restriction Theorems.}
\medskip
 
Let  $F\ss J^2(Z)$ be a subequation on a manifold $Z$, and suppose $i:X\ss Z$ is a 
submanifold.  Then there is a natural induced subequation $i^*F$ on $X$ given by
restriction of 2-jets.  For functions $u\in C^2(Z)$ one has directly that 
\medskip
\centerline
{
$u$ is $F$-subharmonic on Z\qquad $\Rightarrow$\qquad 
$u\bigr|_X$ is $i^*F$-subharmonic on $X$.
}
\medskip
\noindent
Generically this induced subequation $i^*F$ is trivial, i.e., all of $J^2(X)$.  
The first problem is to determine the class of submanifolds for which the restriction
is interesting.   In such cases
we then have the following
\medskip
\noindent
{\bf Question:} When does the implication above hold for all $u\in \USC(Z)$?
\medskip

\noindent
{\bf Example.}  The situation is illustrated by the basic  subequation $\cp$ in $\rn$
whose subharmonics are the convex functions. The restriction of a smooth convex function
$u\in C^\infty(\rn)$ to the unit circle in $\bbr^2$  obeys no proper subequation, while 
the restriction of $u$ to a {\sl minimal} submanifold $M\ss \rn$, of any dimension,
is subharmonic for the Laplace-Beltrami operator on $M$. 
This  assertion carries over to  general convex functions $u$.


\vfill\eject
\noindent
{\bf \II.1.  The First General Theorem.}  The paper [HL$_{9}$]  establishes two  restriction theorems of a general nature, each of which has interesting applications. The first entails the following technical
hypothesis.  Fix  coordinates $z=(x,y)$ on   $Z$ so that locally $X \cong \{y=y_0\}$.

\medskip
\noindent
{\bf The Restriction Hypothesis:}  
Given $x_0\in X$ and  $(r_0, p_0, A_0)\in \Jtn$ and given
$z_\e = (x_\e, \, y_\e)$ and $r_\e$ for a sequence of real numbers $\e$ converging to 0: \medskip
$$
{\rm If\ \ } \left( r_\e,\ \left(p_0+A_0(x_\e-x_0), \ {{y_\e-y_0}\over\e}\right), \ \left(\matrix{A_0&0\cr 0&{1\over \e}I}\right)\right)\ \ \in\ \ F_{z_\e} 
$$
$$
{\rm and\ \ }  x_\e\ \to\  x_0,\ \ {{|y_\e-y_0|^2}\over\e}\ \to\ 0, \ \  r_\e\ \to\ r_0,
$$
then
$$
(r_0, p_0, A_0)\ \in \ (i^*F)_{x_0}.
$$

\smallskip
\noindent
{\bf   Theorem \II.1.1.}  {\sl  Suppose $u\in \USC(Z)$.  Assume the restriction hypothesis and  suppose that
$(i^*F)$ is closed.  Then}
 $$
 u\in F(Z)  \qquad \Rightarrow\qquad u\bigr|_X \in (i^*F)(X).
 $$

If    $(i^*F)$ is not closed, the conclusion holds with 
$(i^*F)$ is replaced by $\overline{(i^*F)}$.


\vskip .3in
\noindent
{\bf \II.2.  Applications of the First General Theorem.}  Theorem \II.1.1 applies to several interesting cases. In the following,
the term {\sl restriction holds} refers to the conclusion of Theorem \II.1.1. The reader is
referred to   [HL$_{9}$] for full statements and proofs.  

\Theorem{\II.2.1} {\sl  Let  $\bbf$ be  a constant coefficient subequation in $\rn$.
Then restriction holds for all affine subspaces $X$ for which $i^*\bbf$ is closed.}

\medskip
More generally, if $u$ is $\bbf$-subharmonic, then $u\bigr|_X$ is $\overline{i^* \bbf}$-subharmonic.

Consider now a  second-order linear operator $\bll$ with smooth coefficients
on $\rn$. Fix linear coordinates $z=(x,y)$ and suppose $X\cong \{y=y_0\}$ as above.
Using the summation convention, write
$$
\bll(u) \ =\ A_{ij}(z) u_{x_i x_j} +a_i(z) u_{x_i} +  \a(z) u   +
B_{k\ell}(z) u_{y_k y_{\ell}} + b_k(z) u_{y_k} + C_{ik}(z) u_{x_i y_k}
$$
Suppose the subequation $L$ corresponding to $\bll u\geq0$ satisfies positivity.
If any one of the coefficients $B(x_0, y_0), b(x_0, y_0)$
or $C(x_0, y_0)$ is non-zero, restriction is trivial locally since
$i^*L$ is everything for $x$ near $x_0$.  Hence, we assume the following
$$
B(x, y_0), b(x, y_0),\ {\rm and}\ C(x, y_0)\ \ {\rm vanish\ identically\ on\ } X
\eqno{(\II.2.1)}
$$

\Theorem{\II.2.1} {\sl  Assuming (\II.2.1), restriction holds for the linear operator
$L$ to $X$.}

\medskip

This result for linear operators proves to be quite useful.

The next result concerns geometric subequations (see Example \BB.1.11) on
general riemannian manifolds $Z$.

\Theorem {\II.2.3} {\sl  Let $\GG \ss G(p,TZ)$ be a closed subset of the bundle of
tangent $p$-planes on $Z$, which admits a fibre-wise neighborhood  retract
(a sub-bundle for example). Let $F(\GG)$ be the induced subequation on $Z$, defined
  as in (\BB.1.6) using the riemannian hessian.  Then restriction holds for all minimal
$\GG$-submanifolds $X\ss Z$, i.e., minimal submanifolds with $T_xX \in \GG_x$ for all
$x\in X$.
}


\vskip .3in
\noindent
{\bf \II.3.  The Second General Theorem.} Let $F$ be a subequation on a manifold $Z$
and fix a submanifold $i:X\ss Z$ as above. In  \CC.2.3 we defined  the notion
of $F$ being locally jet-equivalent to a constant coefficient subequation $\bbf$.  
In our current situation there is a notion of  $F$ being locally jet-equivalent to   $\bbf$
{\sl relative to the submanifold $X$}. This entails $i^*F$ being locally jet-equivalent
to a constant coefficient subequation (assumed closed) on $X$. For details, see [HL$_{9}$, \S\S 9 and 10].

\Theorem{\II.3.1} {\sl
If $F$ is locally jet-equivalent to a constant coefficient subequation relative to 
$X$, then restriction holds for $F$ to $X$.}


\vskip .3in
\noindent
{\bf \II.4.  Applications of the Second General Theorem.}  A nice application of Theorem \II.3.1 is the following.

\Theorem {\II.4.1} {\sl
Let $Z$ be a riemannian manifold of dimension $n$ and $F\ss J^2(Z)$ a subequation canonically
determined by an O$(n)$-invariant constant coefficient subequation $\bbf\ss\bbj^2$.
Then restriction holds for $F$ to any totally geodesic submanifold $X\ss Z$.
}

\medskip

Suppose now that $Z$ is a riemannian manifold with a topological $G$-structure
and $F\ss J^2(Z)$ is  determined by a $G$-invariant 
constant coefficient subequation as in Section \HH.2.
The local framings $e_\a$ appearing in Definition \HH.2.1 are called {\sl admissible}.
So also is any framing of the form $e_\a' = g e_\a$ for a smooth map  $g:U_\a\cap U_\b \to G$.
 A submanifold $X\ss Z$ is said to be {\bf compatible} with the $G$-structure
  if at every point $z\in X$ there is an admissible 
 framing $e$ on a neighborhood $U$ of $z$ such that on  $X\cap U$
 $$
e_1,...,e_n \ \ {\rm are \ tangent\ to\ \ } X\cap U \and e_{n+1},...,e_N \ \ {\rm are \ normal\ to\ \ } X\cap U.
$$
For example, if  $G={\rm U}(N/2)$, then any submanifold of constant CR-rank is compatible.

\Theorem{\II.4.2}  {\sl
Let $Z$ be a riemannian manifold with topological $G$-structure, 
and $F\ss J^2(Z)$ a subequation canonically
determined by a  $G$-invariant constant coefficient subequation $\bbf\ss\bbj^2$.
Then restriction holds for $F$ to any totally geodesic submanifold $X\ss Z$ which is
compatible with the $G$-structure.
}

\medskip
There is a  further application of Theorem \II.3.1  to almost complex manifolds, which is discussed in \S \KK.

\vskip .3in


\noindent{\headfont \JJ.\  Convex Subequations and the Strong Bellman Principle.}
\medskip
 
 An elementary fact, known to all, is that a closed convex set in in a vector space $V$ is the intersection
 of the closed half-spaces containing it.  Put this into a family and you have a fundamental
 principle, which we call the {\sl Bellman Principle}, for dealing with nonlinear pde's 
 which are convex. Specifically, suppose $F\ss J^2(X)$ is a convex subequation --
  one with the property that every fibre $F_x$ is convex. Then, under mild assumptions,
   $F$ can be written locally  as the intersection of a family of  {\bf linear} subequations.  These are subequations of the form
  $$
Lu \ =\   \bra a { D^2u} + \bra b  {Du} + c  u \ \geq\ \l,
  \eqno{(\JJ.1)}
  $$
 where, from the   Conditions
  (P)  and (N) for $F$, 
 one can deduce that the matrix function $a$ 
 and the scalar function $c$ satisfy
 $$
 a\ \geq\ 0
 \and 
 c\ \leq\ 0.
 \eqno{(\JJ.2)}
  $$
The introduction of these  local  linear equations goes back to Richard Bellman
and his work in dynamic programing.  These equations can be found in many areas of mathematics.
Examples close in spirit to those above appear in work of Bedford-Taylor [BT$_*$] and Krylov [K].

It is obviously a big improvement if all  the linear equations in
(\JJ.1) needed to carve out $F$ can be taken to have 
$$
a\ >\ 0, 
 \eqno{(\JJ.3)}
  $$
for then the machinery of uniformly elliptic linear equations  can be brought to bear.

More specifically: any $F$-subharmonic function $u$ is 
locally  a viscosity subsolution of $Lu\geq\l$.
From this one sees that $u$ is a classical subsolution (see [HL$_{10}$, Thm. A.5]), 
and if $a>0$, the results of [HH] apply to prove that $u$ is $\lloc$.
It can then be shown that $u$ is a distributional subsolution to $Lu\geq\l$,
and the full linear elliptic theory ([Ho$_{2}$] or [G] for example) applies.

This naturally raises the question: What assumptions on $F$ will guarantee that it is cut
out by linear equations with $a>0$? 

This question has two parts. The first concerns only the convex geometry of the fibres
$F_x$ at each point $x$; in other words, the question for a convex  
constant coefficient subequaton $\bbf\ss\bbj^2$. The second only 
 involves the mild  regularity condition
that  a containing half-space for $F_x$ extends locally to a linear 
(variable coefficient) subequation
containing $F$.

These questions have been discussed in [K], and an account has also been
given in [HL$_{13}$], where the answer to the first question is given as follows.  
We say that  a subset $C\ss \Symn$  {\sl depends  on all the variables}
if there is no proper subspace $W\ss \rn$ and subset $C'\ss \Sym(W)$ such that 
$A\in C \iff A\bigr|_W \in C'$.  Then a (constant coefficient) subequation 
$
\bbf \ \ss\ \bbj^2 \ =\ \bbr\times \rn\times \Symn
$
is said to {\sl depend weakly on all the second-order variables} if  for each  $(r,p)\in \bbr\times \rn$,
the fibre   $\bbf_{(r,p)} = \{A\in \Symn : (r,p,A)\in \bbf\}$
depends on all the variables.

\Theorem{\JJ.1} {\sl If $\bbf$ depends weakly on all the second-order variables, then
$\bbf$ can be written as the intersection of a family of half-space subequations 
$\bra a A +\bra b p +c r\geq\l$ with $a>0$.}

\Note{\JJ.2} For subequations which do not depend on all the second order
variables, the conclusions above fail.  Consider the (geometrically 
determined) subequation 

\centerline{$\bbf\cong \{u_{xx}\geq0\}$}

\noindent  in the $(x,y)$-plane.  Any continuous
function $u(y)$ is $\bbf$-subharmonic, in fact, $\bbf$-harmonic, but not in general
$\lloc$.
\medskip

See [HL$_{13}$] for a full discussion of these matters.

\vfill\eject


\noindent{\headfont \KK.\  Applications to Almost Complex Manifolds.}
\medskip
 In this section we consider completely general almost complex manifolds
 $(X,J)$ where $J:TX\to TX$ is  smooth bundle map with $J^2\equiv -\Id$.
 On any such manifold there is an {\sl intrinsically defiined} subequation
 $$
 F(J) \ \ss\ J^2(X),
 $$
for which, when the structure is integrable, the $F(J)$-subharmonic functions are
exactly the standard plurisubharmonic functions.  Hence, the results and techniques
discussed in this paper apply to give a full-blown potential theory on almost
complex manifolds, which extends the classical theory.  
The consequences are worked out in detail in [HL$_{10}$].  Here are a few highlights.


\vskip .3in
\noindent
{\bf \KK.1.  $J$-Holomorphic Curves.}  A submanifold $Y\ss X$ is an {\sl almost complex
submanifold} if $J(T_yY)=T_yY$ for all $y\in Y$. In general dimensions such submanifolds
exist only rarely.  However, when the real dimension of $Y$ is two, $Y$
is called a {\bf $J$-holomorphic curve}, and we have the following  important classical result.

\Theorem{\KK.1.1.  (Nijenhuis and Woolf [NW])}
{\sl
For each point $x\in X$ and each complex tangent line $\ell\ss T_xX$, there exists
a $J$-holomorphic curve passing through $x$ with tangent direction $\ell$.
}
\medskip

The restriction result \II.3.1 applies in this case to prove the following.
For historical compatibility we replace the term ``$F(J)$-subharmonic''
with ``$F(J)$-plurisubharmonic''.

\Theorem{\KK.1.2} {\sl
Let $(Y,J_Y)$ be an almost complex submanifold of $(X, J_X)$.  Then the restriction of 
any $F(J_X)$-plurisubharmonic function to $Y$ is $F(J_Y)$-plurisubharmonic.}

\medskip

This leads to the following result equating two natural definitions of plurisubharmonicity.
We recall  that an  almost complex structure $J$ on a 2-dimensional manifold $S$ is
always integrable, and  all notions of (usc) subharmonic functions on $(S,J)$ coincide.  

\Theorem{\KK.1.3} {\sl
A function $u\in \USC(X)$ is $F(J)$-\psh if and only if its restriction to every $J$-holomorphic
curve is subharmonic.
}


\vskip .3in
\noindent
{\bf \KK.2.  Completion of the Pali Conjecture.}  There is a third definition of 
$J$-plurisub-harmonic functions on an almost complex manifold $(X,J)$, which makes
sense for any distribution $u\in \cd'(X)$. Any such distribution $u$ is known to be
$\lloc$.  By work of  Nefton Pali [P] we know that any $\in \USC(X)$ which is
$J$-\psh in the sense of Section \KK.1, is $\lloc$ on $X$ and $J$-\psh as a distribution. 
In the converse direction he showed that if a $J$-\psh distribution 
$u$ has a continuous representative (as a $[-\infty,\infty)$-valued function),
then it is $J$-\psh as above.  He further conjectured that the converse should
hold in general.  This was proved in [HL$_{10}$]. 

The proof used the Strong Bellman Principle and involved showing that the
upper semi-continuous representative of the $\lloc$-class obtained for each
of the associated linear equations, is independent of the linear equation.
It is, in fact, given by the {\sl essential upper-semi-continuous regularization}
$$
u^*_{\rm ess}(x) \ \equiv\ \lim_{r\searrow 0} \left\{{\rm ess} \sup_{B_x(r)} u\right\}
$$
which depends only on the $\lloc$-class of $u$.


\vskip .3in
\noindent
{\bf \KK.3.  The Dirichlet Problem for Maximal Functions.}  Theorem 8.12 applies
 in this case to prove existence and uniqueness  for the Dirichlet problem
for $J$-maximal functions.
One can show that  the more classical notion of  a  function  $u$  being  {\bf $J$-maximal} 
(going back to [B], [W]), is the same as $u$ being
$F(J)$-harmonic, i.e., $u$ is
  $F(J)$-(pluri)subharmonic and $-u$ is $\wt F(J)$-subharmonic. A domain $\O\ss\ss X$
with smooth boundary is strictly {\sl $J$-convex} if it has a strictly $F(J)$-\psh defining function.

\Theorem{\KK.3.1} {\sl
Let  $\O\ss\ss X$ be a strictly $J$-convex domain in an almost complex manifold $(X,J)$.
Then the Dirichlet problem for $J$-maximal functions in uniquely solvable on $\O$ for 
all continuous boundary values $\vf\in C(\bo)$. 
}

\noindent
{\bf Note \KK.4.1.}  Recently Szymon Pli\'s  has also studied the Dirichlet problem on almost complex
manifolds [Pl].  His result  is the almost-complex analogue of a main 
result in [CKNS].  It  treats the  inhomogeneous Monge-Amp\`ere equation
 with positive right hand side. All data are assumed to be smooth,
 and  complete regularity is established for the solution. 
\medskip

\vskip .3in


\noindent{\headfont Appendix A.   A Pocket Dictionary.}
\medskip

  The conventions adopted in this  paper (and related ones)
 are not common in the literature, but they have  advantages, particularly
 for applications to calibrated geometry and to branches of polynomial operators.
 For readers hard-wired to standard notation (as in, say,  [CIL]), we give 
 here  a concise translation of concepts to serve as a guide. 
 
 Classically, a fully nonlinear partial differential equation for a smooth function 
 $u(x)$ on an open set
 $X \ss \rn$ is written in the form 
 $$
 f(x, u,  Du, D^2u)\ =\ 0
 $$
 for a given  contiinuous  function
 $
 f: X \times \bbr\times\rn\times\Symn \ \arr\ \bbr.
 $
 
Here the function $f$  is typically   replaced by the closed set 
 $$
 F\ \equiv\ \{(x,r,p,A) : f(x,r,p,A)\geq 0\}.
 $$
For $C^2$-functions $u(x)$  we have the following translations. Set $J_x^2u \equiv (x,u,Du,D^2u)$.
$$
\eqalign
{
u\ \ {\rm is \  a \ {\bf subsolution}\ } \qquad   &{\bf <\!\!\!- \!\! \!- \!\!\!-\!\!\!>}    \qquad u \ {\rm is\ } F \, {\bf subharmonic}, {\ \rm i.e.,}
\cr
f(x,u,Du,D^2u)\geq0    \qquad   &{\bf <\!\!\!- \!\! \!- \!\!\!-\!\!\!>} \qquad  J_x^2 u\in F \ \ \ \forall x\in X.  
}
$$

$$
\eqalign
{
u\ \ {\rm is \  a \ {\bf supersolution}\ } \qquad   &{\bf <\!\!\!- \!\! \!- \!\!\!-\!\!\!>} \qquad -u \ {\rm is\ } \wt F \, {\bf subharmonic}, {\ \rm i.e.,}
\cr
f(x,u,Du,D^2u)\leq0    \qquad   &{\bf <\!\!\!- \!\! \!- \!\!\!-\!\!\!>} \qquad  -J_x^2 u\in \wt F \ \ \ \forall x\in X.  
}
$$

$$
\eqalign
{
u\ \ {\rm is \  a \ {\bf solution}\ } \qquad   &{\bf <\!\!\!- \!\! \!- \!\!\!-\!\!\!>} \qquad u \ {\rm is\ } F \, {\bf harmonic}, {\ \rm i.e.,}
\cr
f(x,u,Du,D^2u) = 0    \qquad   &{\bf <\!\!\!- \!\! \!- \!\!\!-\!\!\!>} \qquad  J_x^2 u\in \partial F \ \ \ \forall x\in X
\cr
&{\bf <\!\!\!- \!\! \!- \!\!\!-\!\!\!>} \qquad 
 u \ \ {\rm is\ } F\, {\rm subharmonic \  and\  }\cr
&  \qquad\qquad
 -u \ \ {\rm is\ } \wt F\, {\rm subharmonic}
 }
$$
These same translations apply to any upper semi-continuous function $u$ by applying them
to test functions at each point $x$.
\medskip

We also have the following translations between some of the standard structural conditions
placed on the function $f$ and conditions on the set $F$.
Let $\cp \equiv \{(0,0,A) : A\geq0\}$ and $\cn \equiv \{(r,0,0) : r \leq0\}$.

 $$
\eqalign
{
f\ \ {\rm is \  {\bf degenerate \  elliptic}\ } \qquad  \qquad &{\bf <\!\!\!- \!\! \!- \!\!\!-\!\!\!>} \qquad F \ {\rm satisfies\ }  {\bf positivity}, {\ \rm i.e.,}
\cr
f(x,r,p, A+P)\geq f(x,r,p, A) \ \forall P\geq0   \qquad   &{\bf <\!\!\!- \!\! \!- \!\!\!-\!\!\!>} \qquad \qquad F+\cp\ \ss\ F.  
}
$$

 $$
\eqalign
{
f\ \ {\rm is} \  {\bf monotone\ in\ the \ dependent\  variable}  \qquad &{\bf <\!\!\!- \!\! \!- \!\!\!-\!\!\!>} \qquad F \ {\rm satisfies\ }  {\bf negativity}, {\ \rm i.e.,}
\cr
f(x,r-s,p, A)\geq f(x,r,p, A) \ \forall s\geq0   \qquad   &{\bf <\!\!\!- \!\! \!- \!\!\!-\!\!\!>} \qquad \qquad F+\cn\ \ss\ F.  
}
$$

 $$
\eqalign
{
\qquad
f\ \ {\rm is  \  {\bf proper} \ if\ both \ conditions \ hold}  \qquad &{\bf <\!\!\!- \!\! \!- \!\!\!-\!\!\!>} \qquad  
F + \cp \ss F\ \ {\rm and}\ \ F + \cn\ss F
}
$$

 $$
\eqalign
{
f\ \ {\rm is } \  {\bf uniformly\ elliptic }  \qquad &{\bf <\!\!\!- \!\! \!- \!\!\!-\!\!\!>} \qquad  
\cases
{
F + \bbp_{\l,\L} \ \ss\ F \ \ {\rm for\ some\ } 0<\l < \L,  {\rm or\ equivalently,} \cr 
F + \bbp(\d) \ \ss\ F \ \ {\rm for\ some\ } \d>0.
}
}
$$
Here  $\bbp_{\l,\L}$ is the Pucci cone discussed in \S \DD.5, and $\cp(\d)$ is the cone defined in
Example \DD.3.3.
 
 \medskip
 
 It is important to realize that  these translations are not precise equivalences (although 
 there is  an implication).  In passing from the function $f$ to the set $F \equiv \{f\geq0\}$, 
 the behavior of 
 $f$ away from its zero-set is lost. Matters become simpler in a sense.  There are also
 natural examples where the set $\{f\geq0\}$ is not really what one wants to take for the set $F$, and the topological condition required in the ``set'' point of view
 easily corrects matters (see Comment 2 below).
 \medskip

 \medskip
 \noindent
 {\bf Comment 1.} As noted above,  these translations are not equivalences in general.
 For example, the positivity condition $F+\cp\ss F$ is equivalent to the assumption that
 $$
 f(x,r,p,A) \ \geq\ 0
 \qquad\Rightarrow\qquad
 f(x,r,p,A+P) \ \geq\ 0 \ \ \forall \ P\geq0.
 $$
 which is weaker than the inequality on $f$ required for degenerate ellipticity.
 The negativity condition  $F+\cn\ss F$ is equivalent to the assumption that
 $$
 f(x,r,p,A) \ \geq\ 0
 \qquad\Rightarrow\qquad
 f(x,r-s,p,A) \ \geq\ 0 \ \ \forall \ s\geq0.
 $$
 which is weaker than the properness condition placed on $f$ above.
 
  \medskip
 \noindent
 {\bf Comment  2.} The Topological Condition (T) that $F=\overline{\Int F}$, holds for
 most classical equations of interest.  However, there are cases where it fails, such as the 
 infinite Laplaican   $f(p,A) = \bra{Ap}p$ or the $k$-Laplacian $|p|^2 +(k-2) \bra{Ap}p$,
 ($1\leq k\neq 2$).  When it fails,  it is  condition (T)
 that selects the ``correct''  subequation  $F$.

  \medskip
 \noindent
 {\bf Comment  3 (Concerning the first three translations above).} 
 There is an important difference between $u$ being a supersolution and $-u$ being 
 $\ft$-subharmonic, which arises when $\Int F  \neq \{f>0\}$.  Since $ \{f>0\} \ss \Int F$
 (equivalently $\sim \Int F \ss \{f\leq0\}$) we have that 
 $$
 -v \ \ {\rm is \ }  \ft \,{\rm subharmonic} \qquad\Rightarrow\qquad
 v\ \ {\rm is \ an \ } f\,{\rm supersolution}.
 \eqno{(A.1)}
 $$
The fact that the converse is not true is important.  For
 a constant coefficient, pure second-order subequation $F\ss\Symn$,
the more restrictive condition on $v$ in (A.1) ensures that  comparison holds.
That is, with $u$ $F$-subharmonic and $-v$ $\ft$-subharmonic, 
$$
u\ \leq\ v \ \ {\rm on}\ \partial K 
\qquad\Rightarrow \qquad
u\ \leq\ v \ \ {\rm on}\  K 
$$
(See [HL$_{4}$] for a proof.)

 \vskip.3in

\noindent{\headfont Appendix B.   Examples of Basic Monotonicity Cones.}
\medskip

The following is a list of constant-coefficient convex cone subequations   $\bbm$ such that
 the Zero Maximum Principle (see \S \DD.1) holds for $\wt \bbm$-subharmonic functions. 
 In cases (1), (5) and (6) the full maximum principle holds, since these equations are 
 independent of the $r$-variable.

\medskip
\noindent
(1)\ \ $\bbm \ =\ \bbr\times\rn\times \cp$.
Here the $\wt \bbm$-subharmonic functions are the subaffine functions (see Proposition \BB.1.7).
This is a monotonicity subequation for any pure second-order subequation $\bbf= \bbr\times\rn\times \bbf_0$.

\medskip
\noindent
(2)\ \ $\bbm \ =\ \bbr_-\times\rn\times \cp$.
Here  one can characterize the $\wt \bbm$-subharmonics 
as being ``sub''  the functions of the form $\max\{0,a(x)\}$ with $a(x)$ affine (the {\sl affine-plus functions}).
This is a monotonicity subequation for any gradient-independent subequation.

\medskip
\noindent
(3)\ \ $\bbm \ =\ \bbr_-\times \cd \times \cp$ with $\cd \ss \rn$ a ``directional'' convex cone with vertex 
at the origin and non-empty interior.

\medskip
\noindent
(4)\ \ $\bbm\ =\  \{(r,p,A)
\in\bbj^2 :  r\leq - \g |p|, \ \ p\in \cd\ \ {\rm and\ }\  A\geq0\}$
with $\g>0$ and $\cd\ss\rn$ as above.

\medskip
\noindent
(5) \ \ $\bbm\ =\ \bbr\times \bbm_0$\ \  \ with\ \ \  $(p,A)\in \bbm_0\ \iff\ \bra {Ae}e -\l|\bra p e|\geq0 \ \ \forall \, |e|=1$

\medskip

For the next example the Maximum Principle only  holds for  compact sets $K\ss\rn$ which are
contained in a ball of radius $R$.

\medskip
\noindent
(6) \ \ $\bbm\ =\ \bbr\times \bbm_0$\ \  \ with\ \ \  $(p,A)\in \bbm_0\ \iff\   A-{|p|\over R}{\rm Id} \geq0$

\medskip

The proofs depend on the following elementary result.

\Theorem{B.2}  {\sl
Suppose $\bbm$ is a constant coefficient convex subequation and $K\ss\rn$ is compact.
If $K$ admits a smooth function $\psi$ which is strictly $\bbm$-subharmonic on $\Int K$,
then the Zero Maximum Principle holds for the dual subequation $ \wt \bbm$ on $K$.
}

\pf Suppose that the (ZMP) fails for $u\in \USC(K)$.
We will show that  there exists a point $\bar x \in \Int K$ and $\e>0$ such that $\vf \equiv -\e\psi$
is a test function for $u$ at $\bar x$. This proves that $u$ is not $\wt\bbm$-subharmonic
 near $\bar x$ because $J^2_{\bar x}\psi \in \Int \bbm$ implies that $J^2_{\bar x}\vf 
 =- \e J^2_{\bar x}\psi  \notin \wt \bbm$.
 
 By assumption, $u\leq 0$ on $\partial K$ but $\sup_K u>0$.  The negativity condition
 (N) for $\wt \bbm$ allows us to subtract a small number from $u$ and assume that
 $u<0$ on $\partial K$ with  $\sup_K u>0$. Set $v\equiv u+\e \psi$.  
 Then with $\e>0$ sufficiently small, $v<0$ on $\partial K$ but $\sup_K v>0$.
Now let $\bar x$ denote a maximum point for $v$ on $K$.
Since $\bar x \in \Int K$, this proves that $\vf \equiv -\e \psi$ is a test function for $u$
at $\bar x$ as desired. \qed

\medskip
\noindent
{\bf Proof of (1) -- (4)}.
 Since the $\bbm$ in (4) is contained in the other three $\bbm$'s, it suffices to find a 
 strictly $\bbm$-subharmonic function for $\bbm$ defined as in (4).  
 Choose  $\psi(x) \equiv \half \d |x-x_0|^2 -c$ with $\d,c>0$.  Denote the jet coordinates
 of $\psi$ at $x\in K$ by $r=\psi(x)$, $p= \d(x-x_0)$ and $A= \d I$.  Choose $x_0\in\rn$
 so that $K\ss x_0+\Int \cd$.  Then $A\in \Int \cp$, $p\in \Int \cd$ and $r+\g |p| 
 = \half \e |x-x_0|^2 - c + \g \d |x-x_0|< 0$ if $c$ is large.\qed

\medskip
\noindent
{\bf Proof of (5).}  
Consider $\psi(x) \equiv {1\over N+1} |x|^{N+1}$. Then one computes that
$$
p\ =\ D\psi \ =\ |x|^N{x\over |x|}
\and
A\ =\ D^2\psi \ =\ |x|^{N-1}\left( I + (N-1)P_{[x]}    \right)
$$
where $P_{[x]}$ is orthogonal projection onto the $x$-line.
Then with $|e|=1$ we have
$$
{1\over |x|^{N-1}}(\bra {Ae} e - \l|\bra pe | ) \ =\ 1 -\l|x|t  +(N-1)t^2   \ \equiv \ g(t).
$$
with $t\equiv | \bra{ {x\over |x|} } e |$.  
We can assume that $0\notin K$ and $x\in K$ implies $|x|\leq R$.
The quadratic $g(t)$ has a minimum 
at $t_0 = {\l|x| \over 2(N-1)}$ with the minimum value $g(t_0) = 1- {\l^2  |x|^2 \over 4(N-1)}
\geq   1- {\l^2  R^2 \over 4(N-1)}$.  Choose $N$ large enough so that this is $>0$.\qed

\medskip
\noindent
{\bf Proof of (6).}  This is similar to the proof of (5).  It reduces to showing that
$g(t) = 1-{|x|\over R}+(N-1)t^2>0$. Now the minimum value (at $t=0$) is $1-{|x|\over R}$.
For the counterexample, consider 
$$
u(x)\ \equiv\ \cases{-(R-|x|)^3 \quad |x|\leq R\cr
\ \ \qquad 0 \qquad\quad  |x|\geq R
}
$$

\vskip .3in



\centerline{\bf References}

\vskip .2in

\noindent
\item{[A$_1$]}   S. Alesker,  {\sl  Non-commutative linear algebra and  plurisubharmonic functions  of quaternionic variables}, Bull.  Sci.  Math., {\bf 127} (2003), 1-35. also ArXiv:math.CV/0104209.  

\smallskip

\noindent
\item{[A$_2$]}   \ \----------,   {\sl  Quaternionic Monge-Amp\`ere equations}, 
J. Geom. Anal., {\bf 13} (2003),  205-238.

 ArXiv:math.CV/0208805.  

\smallskip

\noindent
\item{[AV]}    S. Alesker and M. Verbitsky,  {\sl  Plurisubharmonic functions  on hypercomplex manifolds and HKT-geometry}, arXiv: math.CV/0510140  Oct.2005

\smallskip


 \smallskip

\noindent
\item{[Al]}     A. D. Alexandrov,  {\sl  The Dirichlet problem for the equation Det$\| z_{i,j}\| = \psi(z_1,...,z_n,x_1,...,x_n)$}, I. Vestnik, Leningrad Univ. {\bf 13} No. 1, (1958), 5-24.

\smallskip

 \smallskip

\noindent
\item{[AGV]} M. E. Amendola, G. Galise and A. Vitolo, 
{\sl Riesz capacity, maximum principle and removable sets of fully
nonlinear second-order elliptic operators},  Preprint, University of Salerno.

\smallskip

\noindent
\item{[AFS]}  D. Azagra, J. Ferrera and B. Sanz, {\sl Viscosity solutions to second order partial differential
equations on riemannian manifolds}, ArXiv:math.AP/0612742v2,  Feb. 2007.
 \smallskip

\noindent
\item{[BT$_1$]}   E. Bedford and B. A. Taylor,  {The Dirichlet problem for a complex Monge-Amp\`ere equation}, 
Inventiones Math.{\bf 37} (1976), no.1, 1-44.

\smallskip

\noindent
\item{[BT$_2$]}  \ \----------,   {Variational properties of the  complex Monge-Amp\`ere equation, I. Dirichlet principle}, 
Duke  Math. J. {\bf 45} (1978), no. 2, 375-403.

\smallskip

\noindent
\item{[BT$_3$]}   \ \----------,   {A new capacity for plurisubharmonic functions}, 
Acta Math.{\bf 149} (1982), no.1-2, 1-40.

\smallskip

 \item{[B]}  H. J. Bremermann,
    {\sl  On a generalized Dirichlet problem for plurisubharmonic functions and pseudo-convex domains},
          Trans. A. M. S.  {\bf 91}  (1959), 246-276.
\medskip

 \item{[BH]}  R. Bryant and F. R. Harvey,
    {\sl  Submanifolds in hyper-K\"ahler geometry},
          J. Amer. Math. Soc. {\bf 1}  (1989),  1-31.
\medskip

\item {[CKNS]} L. Caffarelli, J. J. Kohn,  L. Nirenberg, and J. Spruck, {\sl  The Dirichlet  problem for non-linear second order elliptic equations II: Complex Monge-Amp\`ere and uniformly elliptic equations},    Comm. on Pure and Applied Math. {\bf 38} (1985), 209-252.

\smallskip

 \noindent
\item {[CLN]}  L.  Caffarelli,  Y.Y.  Li and L. Nirenberg {\sl  Some remarks on singular solutions
of nonlinear elliptic equations, III: viscosity solutions,
including parabolic operators}.   ArXiv:1101.2833.

\noindent
 \item{[CNS$_1$]}   L. Caffarelli, L. Nirenberg and J. Spruck,  {\sl
The Dirichlet problem for nonlinear second order elliptic equations. I: Monge-Amp\`ere equation},  
Comm. Pure Appl. Math.  {\bf 37} (1984),   369-402.

 \smallskip

\noindent
 \item{[CNS$_2$]}   \ \----------,  {\sl
The Dirichlet problem for nonlinear second order elliptic equations, III: 
Functions of the eigenvalues of the Hessian},  Acta Math.
  {\bf 155} (1985),   261-301.

 \smallskip

\noindent
 \item{[CNS$_3$]}   \ \----------,   {\sl
The Dirichlet problem for the degenerate Monge-Amp\`ere equation},  Rev. Mat. Iberoamericana
  {\bf 2} (1986),   19-27.

 \smallskip

\noindent
 \item{[CNS$_4$]}  \ \----------,  {\sl  Correction to:
``The Dirichlet problem for nonlinear second order elliptic equations. I: Monge-Amp\`ere equation''},  
Comm. Pure Appl. Math.  {\bf 40} (1987),   659-662.

 \smallskip

\noindent
\item{[CCH]}  A. Chau, J. Chen and W. He,  {\sl  Lagrangian mean curvature flow for entire Lipschitz graphs},  
ArXiv:0902.3300 Feb, 2009.

 \smallskip

\noindent
\item{[CGG$_1$]}  Y.-G. Chen, Y. Giga and S. Goto,  {\sl  Uniqueness and existence of viscosity solutions 
of generalized mean curvature flow equations},  Proc. Japan Acad. Ser. A. Math. Sci  {\bf 65} (1989), 207-210.

 \smallskip

\noindent
\item{[CGG$_2$]}  \ \----------,   {\sl  Uniqueness and existence of viscosity solutions 
of generalized mean curvature flow equations},  J. Diff. Geom. {\bf 33} (1991), 749-789..

 \smallskip

\item{[CY$_1$]}
S.-Y. Cheng and S.-T. Yau,   {\sl On the regularity of the Monge-Ampre 
equation $\det(\partial^2u/\partial x_i\partial x_j)=F(x,u)$}, Comm. Pure Appl. Math. 30 (1977), no. 1, 41Ð68. 

 \smallskip

\item{[CY$_2$]}
  \ \----------,     {\sl The real Monge-Ampre equation and affine flat structures}, 
Proceedings of the 1980 Beijing Symposium on Differential Geometry and Differential Equations, Vol. 1, 2, 3 (Beijing, 1980), 339Ð370, Science Press, Beijing, 1982.
 \smallskip

\noindent
\item{[C]}   M. G. Crandall,  {\sl  Viscosity solutions: a primer},  
pp. 1-43 in ``Viscosity Solutions and Applications''  Ed.'s Dolcetta and Lions, 
SLNM {\bf 1660}, Springer Press, New York, 1997.

 \smallskip

\noindent
\item{[CIL]}   M. G. Crandall, H. Ishii and P. L. Lions {\sl
User's guide to viscosity solutions of second order partial differential equations},  
Bull. Amer. Math. Soc. (N. S.) {\bf 27} (1992), 1-67.

 \smallskip

\noindent
\item{[EGZ]}  P. Eyssidieux, V. Guedj, and A. Zeriahi,  {\sl   Viscosity solutions to degenerate Monge-Amp\`ere
equations},  ArXiv:1007.0076.

 \smallskip

\noindent
\item{[E]}   L. C. Evans,  {\sl   Regularity for fully nonlinear elliptic equations
and motion by mean curvature},  pp. 98-133 
in ``Viscosity Solutions and Applications''  Ed.'s Dolcetta and Lions, 
SLNM {\bf 1660}, Springer Press, New York, 1997.

 \smallskip

\noindent
\item{[ESm]}    L. C. Evans and C. K. Smart,   {\sl   Everywhere differentiability of infinite harmonic functions},  
Berkeley preprint, 2012
 \smallskip

\noindent
\item{[ES$_1$]}   L. C. Evans and J. Spruck,  {\sl   Motion of level sets by mean curvature, I},  
J. Diff. Geom. {\bf 33}  (1991), 635-681.

 \smallskip

\noindent
\item{[ES$_2$]}   \ \----------,  {\sl   Motion of level sets by mean curvature, II},  
Trans. A. M. S.  {\bf 330}  (1992),  321-332.

 \smallskip

\noindent
\item{[ES$_3$]}   \ \----------,  {\sl   Motion of level sets by mean curvature, III},  
J. Geom. Anal.   {\bf 2}  (1992),  121-150.

 \smallskip

\noindent
\item{[ES$_4$]}  \ \----------,   {\sl   Motion of level sets by mean curvature, IV},  
J. Geom. Anal.   {\bf 5}  (1995),   77-114.

 \smallskip

 \noindent
\item {[G]} P. Garabedian, { Partial Differential Equations},    J. Wiley and Sons,  New York, 1964.

\smallskip

\noindent
\item{[G]}   L. G\aa rding, {\sl  An inequality for hyperbolic polynomials},
 J.  Math.  Mech. {\bf 8}   no. 2 (1959),   957-965.

 \smallskip

\noindent
\item{[Gi]}   Y. Giga, {\sl  Surface Evolution Equations -- A level set approach},
Birkh\"auser,  2006.

 \smallskip

\noindent
\item{[Giu]}  
E. Giusti,   Minimal surfaces and functions of bounded variation,  Monographs in Mathematics, 80. BirkhŠuser Verlag, Basel, 1984.
\smallskip

\noindent
\item{[Gr]}  
 M. Gromov,  {\sl Pseudoholomorphic curves in symplectic manifolds},  Invent. Math. 82 (1985), no. 2, 307Ð347.

 \smallskip

\noindent
\item{[H]}    F. R. Harvey,   {\sl Removable singularities and structure theorems for positive currents}. Partial differential equations (Proc. Sympos. Pure Math., Vol. XXIII, Univ. California, Berkeley, Calif., 1971), pp. 129-133. Amer. Math. Soc., Providence, R.I., 1973.

\smallskip

 \noindent 
\item {[HL$_1$]}   F. R. Harvey and H. B. Lawson, Jr,  {\sl Calibrated geometries}, Acta Mathematica 
{\bf 148} (1982), 47-157.

 \smallskip

\item {[HL$_{2}$]}  \ \----------, 
 {\sl  An introduction to potential theory in calibrated geometry}, Amer. J. Math.  {\bf 131} no. 4 (2009), 893-944.  ArXiv:math.0710.3920.

\smallskip

\item {[HL$_{3}$]}  \ \----------,   {\sl  Duality of positive currents and plurisubharmonic functions in calibrated geometry},  Amer. J. Math.    {\bf 131} no. 5 (2009), 1211-1240. ArXiv:math.0710.3921.

\smallskip

\item {[HL$_{4}$]}   \ \----------,  {\sl  Dirichlet duality and the non-linear Dirichlet problem},    Comm. on Pure and Applied Math. {\bf 62} (2009), 396-443. ArXiv:math.0710.3991

\smallskip

\item {[HL$_{5}$]}  \ \----------,    {\sl  Plurisubharmonicity in a general geometric context},  Geometry and Analysis {\bf 1} (2010), 363-401. ArXiv:0804.1316.

\smallskip

\item {[HL$_{6}$]}  \ \----------,   {\sl  Dirichlet duality and the nonlinear Dirichlet problem
on Riemannian manifolds},  J. Diff. Geom. {\bf 88} (2011), 395-482.   ArXiv:0912.5220.
\smallskip

\item {[HL$_{7}$]}  \ \----------,    {\sl  Hyperbolic polynomials and the Dirichlet problem},   ArXiv:0912.5220.
\smallskip

\item {[HL$_{8}$]}  \ \----------,   {\sl  G\aa rding's theory of hyperbolic polynomials},   to appear in  {\sl Communications in Pure and Applied Mathematics}.

\smallskip

\item {[[HL$_{9}$]}  \ \----------, {\sl  The restriction theorem for fully nonlinear subequations}, 
   Ann. Inst.  Fourier (to appear). 
ArXiv:1101.4850.

\smallskip

\item {[HL$_{10}$]}  \ \----------,   {\sl  Potential Theory on almost complex manifolds}, 
Ann. Inst. Fourier (to appear).  ArXiv:1107.2584.
\smallskip

\item {[HL$_{11}$]}   \ \----------,  {\sl  Foundations of  $p$-convexity 
and $p$-plurisubharmonicity in riemannian geometry},  ArXiv:  1111.3895.

\smallskip

\item  {[HL$_{12}$]} \ \----------, {\sl  Removable singularities for nonlinear subequations}.   (Stony Brook Preprint).

\smallskip

\item  {[HL$_{13}$]} \ \----------, {\sl The equivalence of  viscosity and distributional
subsolutions for convex subequations -- the strong Bellman principle},
Bol. Soc. Bras. de Mat. (to appear).  ArXiv:1301.4914.

\smallskip

\item {[HL$_{14}$]}  \ \----------, {\sl  Lagrangian plurisubharmonicity and convexity},  Stony Brook Preprint.

\smallskip

\item {[HL$_{15}$]}  \ \----------, {\sl  Radial subequations, isolated singularities and tangent functions},  Stony Brook Preprint.

\smallskip

\item {[HP$_1$]}   F. R. Harvey and J.  Polking,  {\sl Removable singularities of solutions of linear partial differential equations},  Acta Math. {\bf 125} (1970), 39 - 56.

\noindent
\item{[HP$_2$]} \ \----------,  {\sl Extending analytic objects},
Comm.  Pure Appl. Math. {\bf 28} (1975), 701-727.

 \smallskip

\item {[HH]}  M. Herv\'e and R.M. Herv\'e., {\sl  Les fonctions surharmoniques  dans l'axiomatique de
M. Brelot associ\'ees \`a un op\'erateur elliptique d\'eg\'en\'er\'e},    Annals de l'institut Fourier,  {\bf 22}, no. 2 (1972), 131-145.

\smallskip

 \noindent
\item{[Ho$_1$]}
L. H\"ormander,  An introduction to complex analysis in several variables,  Third edition. North-Holland Mathematical Library, 7. North-Holland Publishing Co., Amsterdam, 1990. 
 
  \smallskip

 \noindent
\item{[Ho$_2$]}
 \ \----------,    The analysis of linear partial differential operators. III. Pseudodifferential operators. Grundlehren der Mathematischen Wissenschaften [Fundamental Principles of Mathematical Sciences], 274. Springer-Verlag, Berlin, 1985. 
  \smallskip

   \noindent
\item{[HM]}    L. R. Hunt and J. J. Murray,    {\sl  $q$-plurisubharmonic functions 
and a generalized Dirichlet problem},    Michigan Math. J.,
 {\bf  25}  (1978),  299-316. 

\smallskip

   \noindent
\item{[I]}    H. Ishii,    {\sl  On uniqueness and existence of viscosity solutions of fully nonlinear second-order elliptic pde's},    Comm. Pure and App. Math. {\bf 42} (1989), 14-45.

\smallskip


   \noindent
\item{[IL]}    H. Ishii and P. L. Lions,    {\sl   Viscosity solutions of fully nonlinear second-order
elliptic partial differential  equations},    J. Diff. Eq.  {\bf 83} (1990), 26-78.

\smallskip

   \noindent
\item{[J$_1$]}    R. Jensen,    {\sl  Uniqueness criteria for viscosity solutions of fully nonlinear 
elliptic partial differential  equations},    Indiana Univ. Math. J. {\bf 38}  (1989),   629-667.

\smallskip

   \noindent
\item{[J$_2$]}  
\ \----------,   {\sl Uniqueness of Lipschitz extensions: minimizing the sup norm of the
gradient} , Arch. Rational Mech. Analysis 123 (1993), 51{74.

\smallskip

  \noindent
\item{[JLM]}  
P. Juutinen, P. Lindqvist,  and J.  Manfredi,  {\sl On the equivalence of viscosity solutions and weak solutions for a quasi-linear equation},  SIAM J. Math. Anal. 33 (2001), no. 3, 699Ð717.

   \smallskip

   \noindent
\item{[K]}    N. V. Krylov,    {\sl  On the general notion of fully nonlinear second-order elliptic equations},    Trans. Amer. Math. Soc. (3)
 {\bf  347}  (1979), 30-34.

\smallskip

\item{[La$_1$]}  D.Labutin,
{\sl Isolated singularities for fully nonlinear elliptic equations},  J. Differential Equations  {\bf 177} (2001), No. 1, 49-76.

 \smallskip
 
\item{[La$_2$]}  D.Labutin, {\sl Singularities of viscosity solutions of fully nonlinear elliptic equations}, 
Viscosity Solutions of Differential Equations and Related Topics, Ishii ed., RIMS K\^oky\^uroku
No. 1287, Kyoto University, Kyoto (2002), 45-57

\smallskip

\item{[La$_3$]}   \ \----------, {\sl Potential estimates for a class of fully nonlinear elliptic equations}, 
Duke Math. J. {\bf 111} No. 1 (2002), 1-49.

\smallskip

\item {[L]}   N. S.  Landkof,   {Foundations of Modern Potential Theory},  Springer-Verlag, New York, 1972.

\smallskip

\item {[Le]}  
P. Lelong, 
Fonctions plurisousharmoniques et formes diffŽrentielles positives,  Gordon and  Breach, Paris-London-New York (Distributed by Dunod Žditeur, Paris) 1968.

\smallskip

   \noindent
\item{[LE]}    Y. Luo and A. Eberhard,    {An application of $C^{1,1}$ approximation to
comparison principles for viscosity solutions of curvatures equations},   Nonlinear Analysis
 {\bf 64} (2006), 1236-1254.

\smallskip

\item{[NTV]}  N. Nadirashvili,
  V. Tkachev and S. Vlùadutü  {\sl
 Non-classical Solution to Hessian Equation from
Cartan Isoparametric Cubic}, ArXiv:1111. 0329.

\smallskip
 
\item {[NW]}  A. Nijenhuis and W. Woolf, {\sl  Some integration problems in almost -complex and complex manifolds}, 
Ann. of Math.,  {\bf 77} (1963), 424-489.

\smallskip

\item {[O]}  A. Oberman, {\sl The  convex envelope is the solution of a nonlinear obstacle problem}, 
Proc. A.M.S. {\bf 135}  (2007), no. 6,  1689-1694.

\smallskip

\item {[OS]}  A. Oberman and L. Silvestre, {\sl The Dirichlet problem for the convex envelope}, 
Trans. A.M.S. {\bf 363}  (2011), no. 11,  5871-5886.

\smallskip

\item {[P]}  N. Pali, {\sl Fonctions plurisousharmoniques et courants positifs de type (1,1)
sur une vari\'et\'e presque complexe}, 
Manuscripta Math.  {\bf 118} (2005), no. 3, 311-337.

\smallskip

\noindent
\item{[PZ]}   S. Peng and D. Zhou, 
{\sl Maximum principle for viscosity solutions on riemannian manifolds},    
ArXiv:0806.4768, June 2008.

\smallskip

\item {[Pl]}  S. Pli\'s, {\sl The Monge-Amp\`ere equation on almost complex manifolds}, 
ArXiv:1106.3356, June, 2011.

\smallskip

\item {[Po$_1$]}  
A. V. Pogorelov,  {\sl On the regularity of generalized solutions of the equation} 
det $(\partial^2 u/\partial x_i \partial x_j) = \phi(x_1, ... , x_n)>0$, Dokl. Akad. Nauk SSSR 200, 1971, pp. 534Ð537. 

\smallskip

\item {[Po$_2$]}  
  \ \----------,     {\sl The Dirichlet problem for the n-dimensional analogue of the Monge-Ampre equation}, Dokl. Akad. Nauk SSSR 201, 1971, pp. 790Ð793. 

\smallskip

\item {[RT]} J. B. Rauch and B. A. Taylor, {\sl  The Dirichlet problem for the 
multidimensional Monge-Amp\`ere equation},
Rocky Mountain J. Math {\bf 7}    (1977), 345-364.

\smallskip

\item {[Sh$_1$]} J.-P. Sha, {\sl  $p$-convex riemannian manifolds},
Invent.  Math.  {\bf 83} (1986), 437-447.

\smallskip

\item {[Sh$_2$]}   \ \----------,    {\sl  Handlebodies and $p$-convexity},
J. Diff. Geom.  {\bf 25} (1987), 353-361.

\smallskip

\item {[\SHF]}  
B. Shiffman,  {\sl Extension of positive line bundles and meromorphic maps.},  Invent. Math. {\bf 15} (1972), no. 4, 332-347.

\smallskip

\item {[S$_1$]}  Z. Slodkowski, {\sl  The Bremermann-Dirichlet problem for $q$-plurisubharmonic functions},
Ann. Scuola Norm. Sup. Pisa Cl. Sci. (4)  {\bf 11}    (1984),  303-326.

\smallskip

\item {[S$_2$]}    \ \----------,    {\sl  Pseudoconvex classes of functions. I.
Pseudoconcave and pseudoconvex sets},
Pacific J. of Math.     {\bf 134}  no. 2  (1988),  343-376.

\smallskip

\item {[S$_3$]}    \ \----------,    {\sl  Pseudoconvex classes of functions. II.
Affine pseudoconvex classes on $\bbr^N$},
Pacific J. of Math.     {\bf 141}  no. 1  (1990),  125-163.

\smallskip

\item {[S$_4$]}    \ \----------,    {\sl  Pseudoconvex classes of functions. III.
Characterization of dual pseudoconvex classes on complex homogeneous spaces},
Trans. A. M. S.     {\bf 309}  no.1  (1988),  165-189.

\smallskip

\item {[S$_5$]}    \ \----------,   {\sl 
Complex interpolation of normed and quasinormed spaces in several dimensions},  I. Trans. Amer. Math. Soc. 308 (1988), no. 2, 685Ð711.

\smallskip

\item {[SYZ$_{1}$]} 
A. Strominger, S.-T.  Yau and E. Zaslow, {\sl Mirror symmetry is T-duality},  Winter School on Mirror Symmetry, Vector Bundles and Lagrangian Submanifolds (Cambridge, MA, 1999), 333Ð347, AMS/IP Stud. Adv. Math., 23, Amer. Math. Soc., Providence, RI, 2001.
\smallskip

\item {[SYZ$_{2}$]} 
  \ \----------,    {\sl Mirror symmetry is T-duality},  Nuclear Phys. B 479 (1996), no. 1-2, 243Ð259.

\smallskip

\item {[So]}  P. Soravia, {\sl  On nonlinear convolution 
and uniqueness of viscosity solutions},
Analysis  {\bf 20}    (2000),  373-386.

\smallskip

\item {[TU]} N. S. Trudinger and J.n I. E. Urbas, {\sl Second derivative estimates for equations of Monge-Amp\`ere type},
Bull. Austral. Math. Soc. {\bf 30}  (1984), 321-334.
\smallskip

\item {[W]}   J.  B. Walsh,  {\sl Continuity of envelopes of plurisubharmonic functions},
 J. Math. Mech. 
{\bf 18}  (1968-69),   143-148.

 \smallskip

\item {[WY]}   D. Wang and Y. Yuan,  {\sl Hessian estimates for special Lagrangian
equation with critical and supercritical phases in general dimensions},
 ArXiv:1110.1417.

 \smallskip

\item {[Wu]}   H. Wu,  {\sl  Manifolds of partially positive curvature},
Indiana Univ. Math. J. {\bf 36} No. 3 (1987), 525-548.
 
  \smallskip

\item {[Yau]} S.-T. Yau, {\sl
On the Ricci curvature of a compact KŠhler manifold and the complex Monge-Ampre equation. I}, 
Comm. Pure Appl. Math. 31 (1978), no. 3, 339Ð411. 

\smallskip

\item {[Y]}  Yu Yuan,  {\sl A priori estimates for solutions of fully nonlinear special lagrangian equations},
 Ann Inst. Henri Poincar\'e  
{\bf 18}  (2001),   261-270.

 \end